\documentclass[aos,MSNbibl,nameyear,seceqn,dvips]{arximspdf}
\usepackage{mathrsfs}
\usepackage{dcolumn}

\doi{10.1214/13-AOS1139} %kopijuoti is PTS
\volume{41}
\issue{4}
\pubyear{2013}
\firstpage{2123}
\lastpage{2148}

\makeatletter

\newcolumntype{d}[1]{D{.}{.}{#1}}

\newproclaim{example}{Example}

\newcommand{\cal}{\mathcal}

\newcommand{\xrightarrow}[1]{\stackrel{#1}{\rightarrow}}

\newtheorem{tm}{Theorem}
\newtheorem{cy}{Corollary}
\newtheorem{pn}{Proposition}

\newcommand{\bbeta}{\bolds{\beta}}

\newcommand{\bbetaa}{\bolds{\beta}}

\newcommand{\bmu}{\bolds{\mu}}
\newcommand{\bmus}{\bolds{\mu}}
\newcommand{\bvare}{\bolds{\varepsilon}}
\newcommand{\bSigma}{\bolds{\Sigma}}
\newcommand{\bzero}{\mathbf{0}}
\newcommand{\bg}{\mathbf{g}}
\newcommand{\bZ}{\mathbf{Z}}
\newcommand{\bM}{\mathbf{M}}
\newcommand{\bY}{\mathbf{Y}}
\newcommand{\bW}{\mathbf{W}}
\newcommand{\bX}{\mathbf{X}}
\newcommand{\bI}{\mathbf{I}}
\newcommand{\bu}{\mathbf{u}}
\newcommand{\bR}{\mathbf{R}}
\newcommand{\bv}{\mathbf{v}}

\makeatother

\begin{document}
\begin{frontmatter}

\title{Marginal empirical likelihood and sure independence feature screening}
\runtitle{Marginal EL and SIS}

\begin{aug}
\author[A]{\fnms{Jinyuan} \snm{Chang}\thanksref{t1}\ead[label=e1]{changjinyuan1986@pku.edu.cn}},
\author[B]{\fnms{Cheng Yong} \snm{Tang}\ead[label=e2]{CHENGYONG.TANG@ucdenver.edu}}
\and
\author[C]{\fnms{Yichao} \snm{Wu}\corref{}\thanksref{t2}\ead[label=e3]{wu@stat.ncsu.edu}}
\runauthor{J. Chang, C. Y. Tang and Y. Wu}
\affiliation{Peking University, University of Colorado Denver
and
North~Carolina~State~University}
\address[A]{J. Chang\\
Guanghua School of Management\\
Peking University\\
Beijing, 100871\\
China\\
\printead{e1}}
\address[B]{C. Y. Tang\\
Business School\\
University of Colorado Denver\\
Campus Box 165\\
PO Box 173364\\
Denver, Colorado 80217-3364\\
USA\\
\printead{e2}}
\address[C]{Y. Wu\\
Department of Statistics\\
North Carolina State University\\
2311 Stinson Drive\\
Raleigh, North Carolina 27695-8203\\
USA\\
\printead{e3}} %adresu isvedimo komanda gale!
\end{aug}

\thankstext{t1}{Supported in part by a National Natural Science Foundation
of China Key Program Grant (No. 11131002) and Center for Statistical
Science at Peking University. This work was completed during Jinyuan
Chang's Ph.D. study at Peking University.}
\thankstext{t2}{Supported in part by NSF Grants DMS-09-05561 and
DMS-10-55210 and NIH/NCI Grant R01-CA149569.}

% HISTORY:
\received{\smonth{10} \syear{2012}}
\revised{\smonth{6} \syear{2013}}

% ABSTRACT
%
\begin{abstract}
We study a marginal empirical likelihood approach in scenarios when the
number of variables grows exponentially with the sample size.
%We systematically quantify the properties of the marginal empirical
%likelihood ratio that can be evaluated at arbitrary values of the
%parameter of interest in its support,
%and
The marginal empirical likelihood ratios as functions of the parameters
of interest are systematically examined, and we find that the marginal
empirical likelihood ratio evaluated at zero can be used to
differentiate whether an explanatory variable is contributing to a
response variable or not. Based on this finding, we propose a unified
feature screening procedure for linear models and the generalized
linear models. Different from most existing feature screening
approaches that rely on the magnitudes of some marginal estimators to
identify true signals, the proposed screening approach is capable of
further incorporating the level of uncertainties of such estimators.
Such a merit inherits the self-studentization property of the
empirical likelihood approach, and extends the insights of existing
feature screening methods.
%The proposed screening approach is adaptive and
%lexible, and extends the insights of existing approaches by
%incorporating both the magnitude and level of uncertainties of a
%marginal estimator.
Moreover, we show that our screening approach is less restrictive to
distributional assumptions, and can be conveniently adapted to be
applied in a broad range of scenarios such as models specified using
general moment conditions. Our theoretical results and extensive
numerical examples by simulations and data analysis demonstrate the
merits of the marginal empirical likelihood approach.\looseness=-1
\end{abstract}

% KEYWORDS

\end{frontmatter}

%s1 #&#
\section{Introduction}
High-dimensional data are more frequently encountered in current
practical problems of finance, biomedical sciences, geological studies
and many more areas. Statistical methods for high-dimensional data
analysis have received increasing interests to deal with
large volume of data containing considerably many features; see
\citet{Buhlmannvan2011}, \citet{Hastieetal2009book} and
\citet{FanLv2009Sinica} for overviews. %, review and discussions.
A fundamental
objective of statistical analysis with high-dimensional data is to
identify relevant features, so that effective models can
be subsequently constructed and applied to solve practical problems.

Recently, independence feature screening methods have been
considered, see, for example, \citet{FanLv2008JRSSB},
\citet{FanSong2010AOS} and \citet{FangFengSongJASA2011} for linear
models, generalized linear models and nonparametric additive models,
respectively. \citet{FanLv2008JRSSB} and \citet{FanSong2010AOS}
performed screening by ranking the absolute values of marginal
estimates of model coefficients, and \citet{FangFengSongJASA2011}
carried out screening by ranking integrated squared marginal
nonparametric curve estimates. \citet{FanSong2010AOS} also
discussed independence screening by examining the magnitudes of the
likelihood ratios. More recently, \citet{Wang2012Bioka} considered
a sure independence screening by a factor profiling approach;
\citet{XueZou2011bioka} studied sure independence screening and
sparse signal recovery; see also \citet{Zhuetal2011JASA} and
\citet{Lietal2012JASA} for recent development using model-free
approaches for feature screening, \citet{Lietal2012AOS} for a
robust rank correcation based approach, and \citet{ZhaoLi2012} for
an estimating equation based feature screening approach.

%Consider generic observation $\{Y_i, \bX_i\}_{i=1}^n$ for scalar
%$Y_i$ and vector $\bX_i=(X_{1i},\ldots, X_{ip})^T\in{\mathbb R}^p$.
%independent, while \citet{FanSong2010AOS} and
%condition for $\bX_i$.

The empirical likelihood approach [Owen (\citeyear{Owen1988,Owen2001})]
is demonstrated
effective in scenarios with less restrictive distributional
assumptions for statistical inferences; see \citet{QinLawless1994AS},
\citet{NeweySmith2004Econca} and reference therein. We refer to
\citet{ChenVan2009test} as a review and discussion of recent
development in the empirical likelihood approach.
The scope of the empirical likelihood approach recently has also been
extended to deal with high-dimensional data; see \citet
{Hjortetal2008AS}, \citet{ChenPengQin2008}, \citet
{TangLeng2010Bioka}, %
\citet{LengTang2010}, and \citet{Changetal2012}. Though demonstrated
effective in statistical inferences, the empirical likelihood approach
encounters substantial difficulty when data dimensionality is high.
More specifically, the data dimensionality $p$ cannot exceed the sample
size $n$ in the conventional empirical likelihood
construction. In addition, $p$ can be at most $o(n^{{1/2}})$
or even slower under which asymptotic properties are established
[\citet
{Hjortetal2008AS,ChenPengQin2008,TangLeng2010Bioka,LengTang2010,Changetal2012}].
Therefore, to practically more effectively apply the empirical
likelihood approach, a pre-screening procedure is necessary to
reduce the candidates of target features.\looseness=-1

In this study, we systematically examine the properties of a
marginal empirical likelihood approach where the available
features are assessed one at a time individually. The marginal
empirical likelihood approach only involves univariate optimizations,
so that it provides a convenient device for both theoretical analysis\vadjust{\goodbreak}
and practical implementation.
%so that we are able to explore and justify that the high data
%dimensionality is not an issue for both theoretical analysis and
%practically implementations.
Our analysis reveals the probabilistic
behavior of the marginal empirical likelihood ratios as functions of
the parameters of interest that can be evaluated at
arbitrary values, which itself is a problem of individual interest
because existing studies of the empirical likelihood approach
generally focus on its properties when evaluated at the truth, or at
values in a small
neighborhood of the truth.
Based on our finding, we propose to conduct feature screening by using
the marginal
empirical likelihood ratio evaluated at zero.
We find that a unified screening procedure can be applied
in both linear
models and generalized linear models.
We also demonstrate how the marginal empirical likelihood approach can
be conveniently adapted to solve a broad range of problems
for models specified by general moment conditions.
Hence, the marginal empirical likelihood approach provides a general
and adaptive procedure for solving a broad class of practical problems
for feature screening.
%An example of applying our approach in
%a framework of longitudinal data analysis by using over-identified
%estimating functions \citet{Quetal2000Bioka} demonstrates the
%advantage
%of the marginal empirical likelihood approach.
Our theoretical analyses show that the proposed screening
procedure based on the marginal empirical likelihood approach is
selection consistent---that is, being able to identify the features
that contribute to the response variable when the number of
explanatory variables $p$ grows exponentially with sample size $n$.

% % Theoretical results and a practical concrete example in
%longitudinal data analysis illustrate the a

%By
%demonstrating how such as models specified by general moment
%conditions, we also illustrate an advantage of our approach %We
%establish the properties of such a procedure using the marginal
%empirical likelihood approach

%Moreover, we also explore and show that the proposed screening
%procedure can be adapted to
%a broader framework to conveniently incorporate general available data
%information.

Our study contributes to the sure independence feature screening for
high-dimensional data analysis from the following two substantial
aspects. First of all, a fundamental difference of our approach to
all existing approaches is that
the marginal empirical likelihood ratio statistic is a self-studentized quantity
[\citet{Owen2001}] while other
existing screening methods generally rely on the ranking of features
based on magnitudes of some marginal estimators.
Therefore, our approach is able to incorporate additionally
%without considering
the level of uncertainties associated with the estimators
to conduct feature screening.
%
%is able to incorporate
This
clearly extends the scope of existing feature
screening approaches by considering more aspects of marginal
statistical approaches.
We show in our simulation studies that when heterogeneity exists in the
conditional variance,
our approach performs much better than a least-squares based approach.
Second, our screening procedure inherits the nonparametric merits of
the empirical likelihood approach.
Specifically, our approach requires no strict distributional
assumptions such as normally distributed errors
in the linear models, or exponential family distributed response in the
generalized linear models. This generalizes the scope and applicability
of our approach. As a result,
we show that the marginal empirical likelihood approach provides a
unified framework for feature screening in linear regression models and
generalized linear models, and
can be conveniently applied for solving a broad class of general problems.

The rest of this paper is organized as follows. We elaborate the
method of the marginal empirical likelihood approach in Section~\ref
{s2}. Properties of the proposed approach are given in Section~\ref
{s3}. Section~\ref{s4} extends the marginal empirical likelihood
approach to a broad framework including models specified by general
moment conditions, and presents an iterative sure screening procedure using\vadjust{\goodbreak}
profile empirical likelihood.
Numerical examples are given in Section~\ref{s5}. We conclude with
some discussions in Section~\ref{s6}. All technical details are
contained in the supplementary material of this paper [\citet
{Changetal2013}].

%s2 #&#
\section{Methodology} \label{s2}

%s2.1 #&#
\subsection{Marginal empirical likelihood for linear models}

Let us motivate the marginal empirical likelihood approach by first
considering the multiple linear regression model
%
%e2.1 #&#
%
\begin{equation}\label{eqlrm}
Y=\bX^{\mathrm{T}}\bbeta+\varepsilon,
\end{equation}
where $\bX=(X_1,\ldots,X_p)^{\mathrm{T}}$ is the vector of explanatory
variables, $\varepsilon$ is the random error with zero mean, and
$\bbeta=(\beta_{1},\ldots,\beta_{p})^{\mathrm{T}}$ is the vector
of unknown
parameters. Hereinafter, we also use $\bbeta$ to denote the truth of
the parameter whenever no confusion arises. Without loss of
generality, we assume hereinafter that the explanatory variables are
standardized such that $\mathbb{E}(X_j)=0$ and
$\mathbb{E}(X_j^2)=1$ $(j=1,\ldots,p)$. For effective and
interpretable practical applications, one may reasonably expect
that among the large number of explanatory variables, only a small
fraction of them contribute to the response variable. We therefore
denote by $\mathcal{M}_*=\{1\leq j\leq
p\dvtx\beta_{j}\neq0\}$ the collection of the effective explanatory
variables in the true sparse model %with
whose size is characterized by its cardinality $s=|\mathcal{M}_*|$.
Here we assume that $s$ is much smaller than $p$, reflecting the case
in many practical applications like in finance, biology and clinical studies.

In the recent literature of high-dimensional data analysis, various
marginal approaches have been applied for locating the true model
$\mathcal{M}_*$; see, for example, \citet{FanLv2008JRSSB},
\citet
{FanSong2010AOS} and \citet{FangFengSongJASA2011}. Among those
approaches, a popular
way is to assess the marginal contribution from a given explanatory
variable $X_j$. Commonly applied criteria for
measuring the marginal contribution are the magnitudes of some
marginal estimators
[\citet{FanLv2008JRSSB,FanSong2010AOS,FangFengSongJASA2011}].
Subsequently, the candidate models are chosen
from the top ranked explanatory variables.

To apply a marginal empirical likelihood approach for the linear
regression model (\ref{eqlrm}), let us consider the marginal moment
condition of the least squares estimator:
%
%e2.2 #&#
%
\begin{equation}\label{eqlrmg}
\mathbb{E}\bigl\{X_{j}\bigl(Y-X_j\beta^M_j
\bigr)\bigr\}=0\qquad(j=1,\ldots,p),
\end{equation}
where $\beta_j^M$ is interpreted as the marginal contribution of
covariate $X_j$ to $Y$.
From~(\ref{eqlrmg}), we can see that $\beta_j^M=\mathbb{E}(X_jY)$ is
the covariance between $X_j$ and $Y$ so that $\beta_j^M=0$ is
equivalent to that $Y$ and $X_j$ are marginally uncorrelated. Here
we note the remarkable difference between $\beta_j^M$ and
$\beta_{j}$ where the latter is the truth of the parameter in
(\ref{eqlrm}).
In general, $\beta_j^M\neq\beta_{j}$ unless ${\mathbb E} (X_i X_j)
=0$ for all $i\neq j$.\vadjust{\goodbreak}
In addition to that from $\beta_{j}$ in the model (\ref{eqlrm}),
$\beta_j^M$ also contains aggregated contribution from other
components that may be correlated with $X_j$. Thus, the correlation
level among covariates has significant impact on the
performance of a screening procedure based on (\ref{eqlrmg}); more
discussions on this are given in a later section containing the main
results.

%Therefore, $\beta_j^M$ contains the aggregated covariance between
%$X_j$ and $Y$ as implied by the model (\ref{eqlrm}).

A marginal empirical likelihood for linear models can be constructed as
follows. Note that
$\mathbb{E}(X_j^2)=1$, therefore (\ref{eqlrmg}) is equivalent to
%
%e2.3 #&#
%
\begin{equation}\label{eqn1}
\mathbb{E}\bigl(X_{j}Y-\beta_j^M\bigr)=0.
\end{equation}
Let $\{(\bX_i,Y_i)\}_{i=1}^n$ be collected independent data,
$g_{ij}(\beta)=X_{ij}Y_i-\beta$ $(j=1,\ldots,p)$ and $X_{ij}$ means
the $j$th
component of the $i$th observation $\bX_i$. Based on (\ref{eqn1}),
we define the following marginal empirical likelihood:
%
%e2.4 #&#
%
\begin{equation}
\label{eqEL} \mathrm{EL}_j(\beta)=\sup\Biggl\{\prod
_{i=1}^nw_i\dvtx w_i\geq0, \sum
_{i=1}^nw_i=1, \sum
_{i=1}^nw_ig_{ij}(\beta)=0
\Biggr\}
\end{equation}
for $j=1,\ldots,p$.
For any given $\beta$ in the convex hull of $\{X_{ij}Y_i\}_{i=1}^n$,
the marginal empirical likelihood ratio is defined as
%
%e2.5 #&#
%
\begin{equation}\label{eqellj}
\ell_j(\beta)=-2\log\bigl\{\mathrm{EL}_j(\beta)\bigr
\}-2n\log n=2\sum_{i=1}^n\log\bigl\{1+
\lambda g_{ij}(\beta)\bigr\},
\end{equation}
where $\lambda$ is the Lagrange multiplier satisfying
%
%e2.6 #&#
%
\begin{equation}\label{eqlambda}
0=\sum_{i=1}^n\frac{g_{ij}(\beta)}{1+\lambda
g_{ij}(\beta)}.
\end{equation}

%s2.2 #&#
\subsection{Extended coverage to generalized linear models}

A merit of the\break marginal empirical likelihood approach is that the
formulation by (\ref{eqEL}) and (\ref{eqellj}) only requires the
moment condition (\ref{eqn1}), rather than specific distributional
assumption of $\varepsilon$ in model (\ref{eqlrm}). This entitles
our approach robustness against the violation of distributional model
assumptions, %Furthermore, the marginal empirical likelihood approach
and thus it can be extended
and adapted to a broader framework. Now we elaborate how the above
marginal empirical likelihood approach can be equally applied when the
response variable $Y$ is in the
exponential family with the density function taking the canonical
form [\citet{GLM1989MN}]:
%
%e2.7 #&#
%
\begin{equation}
\label{eqglm}
f(y)=\exp\bigl\{y\theta-b(\theta)+c(y)\bigr\}
\end{equation}
for some suitable known functions $b(\cdot)$, $c(\cdot)$ and
canonical parameter $\theta$. Further extensions of the marginal
empirical likelihood
approach are discussed in a later section. We refer to \citet
{Kolaczyk1994Sinica} and \citet{ChenCui2003Sinica} for conventional
applications of the empirical likelihood to generalized linear
models. Following the convention
of generalized linear models, we denote the mean function by $
\mu=\mathbb{E}(Y|\bX)=b'(\theta)$ where\vadjust{\goodbreak} $\theta$ is modeled by a
linear function $\beta_0+\bX^{\mathrm{T}}\bbeta$ with
$\bbeta=(\beta_1,\ldots,\beta_p)^{\mathrm{T}}$, and use $V(\mu)$
to denote the
variance of $Y$ expressed as a function of $\mu$. %The density
%function (\ref{eqglm}) implies that
%for some function $g(\cdot)$.

For any $j=1,\ldots,p$, the moment condition based on the marginal
likelihood approach in \citet{FanSong2010AOS} for $\beta_j$ is
%
%e2.8 #&#
%
\begin{equation}\label{eqglmg}
\mathbb{E} \biggl\{\frac{Y-\mu_j}{V(\mu_j)}\,\frac{\partial\mu
_j}{\partial\beta_j} \biggr\}=0,
\end{equation}
where $\mu_j=b'(\beta_0+\beta_j^M X_j)$ is the implied mean function
that is modeled marginally only using $X_j$. Here the $\beta_j^M$ is
again interpreted as the marginal contribution of $X_j$ to the
response variable $Y$; see also \citet{FanSong2010AOS}. By the property
of the exponential family
distribution, $\frac{\partial\mu}{\partial\beta_j}=X_j b''(\theta)$
and $V(\mu)=b''(\theta)$. Then (\ref{eqglmg}) becomes
%
%e2.9 #&#
%
\begin{equation}\label{eqglmg2}
\mathbb{E}\bigl\{X_j(Y-\mu_j)\bigr\}=0.
\end{equation}
%
%The marginal version of (\ref{eqglmg2}) is
For linear models, $b(\theta)=\frac{\theta^2}{2}$, then
$b'(\theta)=\theta$, so that (\ref{eqglmg2}) becomes $
\mathbb{E}\{X_j(Y-X_j\beta_j^M)\}=0 $ by noting that $\beta_0=0$ in
the linear model case, which
is exactly (\ref{eqlrmg}). Hence, (\ref{eqglmg2}) is a natural extension
of (\ref{eqlrmg}) in the generalized linear models.

One way to apply the marginal approach can
be generalizing the definition in (\ref{eqEL}) to be $g_{ij}(\beta)=
X_{ij}\{ Y_i- b'(\beta_0+\beta X_{ij})\}$ $(j=1,\ldots, p)$. %, %the
%marginal empirical likelihood for the contribution of $X_j$ in the
%generalized linear model case can also be constructed with the form
%in (\ref{eqEL}).
However, such a modification is actually not necessary.
To see this, we note that when the marginal contribution $\beta
_j^M=0$, then the
marginal moment condition (\ref{eqglmg2}) becomes $\mathbb{E}[X_j\{
Y-b'(\beta_0)\}]=0$.
%mean function of $Y$ becomes $b'(\beta_0)$.
Hence, it implies that the covariance between
$X_j$ and $Y$ is 0, which exactly shares the same implication of (\ref
{eqn1}) as in the linear models.
From this perspective, (\ref{eqglmg2}) and (\ref{eqn1}) are
essentially equivalent.
Additionally, the response variable in practice can always be centered
to have zero mean. This fact eliminates the concern on the intercept
$\beta_0$ in the generalized linear models when considering a marginal
empirical likelihood approach. As a result, %the marginal empirical
%likelihood (\ref{eqEL}) with the same $g_{ij}(\beta)=X_{ij} Y - \beta$
%equally applies in the generalized linear models for centered response
%variable $Y$. %Therefore,
we conclude that a unified marginal empirical likelihood construction
(\ref{eqEL}) with the same $g_{ij}(\beta)=X_{ij} Y_i - \beta$
can be equally applied for both linear models and
generalized linear models with centered response variable $Y$.
The implication of this unified construction is also intuitively very
clear by interpreting $\beta$ as the covariance between a covariate
and the response variable.

Furthermore, we note that the distributional assumption (\ref{eqglm}) is
actually not required in our marginal empirical likelihood
approach. Therefore our approach is not restricted to the exponential
family (\ref{eqglm}). Since we only require the marginal moment
condition (\ref{eqglmg2}), our approach can be applied with the
quasi-likelihood approach and it also works with misspecified
variance functions [\citet{GLM1989MN}].

%Thanks to the merit of being self-studentized
%marginal empirical likelihood approach is able to incorporate both the
%magnitudes and the level of uncertainties associated with the marginal
%approaches. %In addition, the screening procedure only requires the
%marginal moment condition (\ref{eqn1})

%We remark that
The marginal empirical likelihood ratio (\ref{eqellj}) with
$g_{ij}(\beta)=X_{ij}Y_i-\beta$ evaluated at $\beta=0$---that is,
$\ell_j(0)$---has a very clear
practical interpretation by noting that it can be used to test the null
hypothesis $H_0\dvtx\beta_j^M=0$.\vadjust{\goodbreak}
By noting additionally the intuitively clear fact that $\ell_j(0)$
should not be large if $\beta_j^M=0$, we can see
that $\ell_j(0)$ can be used as a device for feature screening.
More specifically, we have the following procedure:

Step 1: Evaluating $\ell_j(0)$ for all $j=1,\ldots,p$, where
$\ell_j(\cdot)$ is defined in (\ref{eqellj}) with
$g_{ij}(\beta)=X_{ij}Y_i-\beta$. If $0$ is not in the convex hull of
$\{X_{ij}Y_i\}_{i=1}^n$, we define $\ell_j(0)=\infty$ as a strong
evidence of significance in predicting $Y$ using $X_j$.

Step 2: Given a threshold level $\gamma_n$, select a set of variables
by
\[
\widehat{\mathcal{M}}_{\gamma_n}=\bigl\{1\leq j\leq p\dvtx\ell_j(0)
\geq\gamma_n\bigr\}.
\]
We specify in the next section the requirement for $\gamma_n$ so
that the screening procedure is consistent. On the other hand, however,
explicitly
identifying $\gamma_n$ in practice is generally difficult because it involves
unknown constants. Thus, a screening procedure can be practically
implemented in a way such that $\widehat{\mathcal{M}}_{\gamma_n}$
recruits candidate features until certain size such as
$n^{{1/2}}$ is achieved.

We remark that the evaluation of $\ell_j(\beta)$ in (\ref{eqellj})
in practice is actually very easy by
noting that all optimizations involved are univariate, which is very
convenient for practical applications.
On the other hand, our procedure only needs to evaluate the marginal
empirical likelihood ratio (\ref{eqellj}) at $\beta=0$ and
avoids the estimation of $\beta_j^M$ when
conducting the feature screening. %
\section{Main results} \label{s3}

Now we present main results for the marginal empirical likelihood
ratio in (\ref{eqellj}) with the unified specification
$g_{ij}(\beta)=X_{ij}Y_i-\beta$ that are generally applicable for
both
linear models and generalized linear models.
In our discussion hereinafter, let $\rho_j={\mathbb E}(X_j Y)$.
If
$\rho_j=0$, it is well known that $\ell_j(0)$ is asymptotically
chi-square distributed with 1 degree of freedom [Owen (\citeyear
{Owen1988,Owen2001})].
If $\rho_j\neq0$, however, the properties of $\ell_j(0)$ is
generally less
clear, which is also a question of independent interest.
Specifically, if $\beta=\rho_j+\tau\sigma n^{-{1}/{2}}$ where
$\sigma^2=\operatorname{var}(X_j Y)$, %and $g_{ij}(\beta)=X_{ij}Y_i-
%%\blue{(what is $g_{ij}$?)},
it can be shown following the same argument of \citet{Owen1988} that
$
\ell_j(\beta)\xrightarrow{d}\chi_1^2(\tau^2) %\textrm{as} n
$ as $n\to\infty$
under some regularity conditions where $\tau^2$ is a noncentrality
parameter. But if $\beta-\rho_j$ converges to zero at a rate
slower than $n^{-{1}/{2}}$, the exact diverging rate of $\ell
_j(\beta)$ is less clear in existing literature.

We first present a general result that shows that the empirical
likelihood ratio
$\ell_j(\beta)$ is no longer $O_p(1)$ when $\beta-\rho_j$ converges
to $0$ but $n^{{1/2}}(\beta-\rho_j)$ diverges.
%This result broadly applies in the generalized linear models.

%
%pr1 #&#
%
\begin{pn}\label{pn1}
Suppose that $U_1,\ldots,U_n$ are independent and identically
distributed random variables with $\mathbb{E}(|U_i|^\nu)<\infty$ for
some $\nu\geq3$. Replacing $g_{ij}(\beta)$ in (\ref{eqellj}) and
(\ref{eqlambda}) by $U_i-\mu$ for all $i=1,\ldots,n$, we obtain
$\ell(\mu)$. If $|\mu-\mu_0|=O(n^{-w})$ for some
$w\in(\frac{1}{\nu},\frac{1}{2})$, then
\[
\frac{\ell(\mu)}{n(\mu-\mu_0)^2\sigma^{-2}}\xrightarrow{p}1\qquad\mbox
{as } n\rightarrow\infty,
\]
where $\mu_0=\mathbb{E}(U_i)$ and
$\sigma^2=\mathbb{E}\{(U_i-\mu_0)^2\}$.\vadjust{\goodbreak}
\end{pn}

We note that \citet{ChenGaoTang2008} contains a related result showing
that the empirical likelihood ratio is diverging when evaluated at
values far enough from the truth. Our Proposition \ref{pn1} contains
the specific diverging rate of the empirical likelihood ratio.
Proposition \ref{pn1} implies that if $\beta- \rho_j$ converges
to zero at a rate slower than $n^{-{1}/{2}}$,
$\ell_j(\beta)=O_p\{n(\beta-\rho_j)^2\}$.
On the other hand, if
$\beta-\rho_j$ does not weaken to zero, our Theorem \ref{tm1} presented later
shows that $\ell_j(\beta)$ has high probability to
take large value.
On the other hand, as clearly shown in our proof of Proposition \ref
{pn1} given in \citet{Changetal2013}, the statistics $\ell_j(0)$ is
self-studentized, and hence it incorporates the level of uncertainties
from using the finite sample moment conditions.
Such a feature is desirable because in practice levels of
uncertainties corresponding to different covariates can be different
when contributing to the response variable of interest. This may
confound the ranking for feature screening based on marginal
estimators themselves without considering their standard errors,
not mentioning incorporating the level of uncertainties is
difficult especially when handling high-dimensional statistical
problems.

% and on the other hand extend the scope of an independence screening
%procedure by automatically incorporating the level of uncertainties.
%%If the
%identification condition proposed as
%$\min_{j\in\mathcal{M}_*}\mathbb{E}X_jY\geq O(n^{-\kappa})$ for some
%$\kappa\in(0,1/2)$, for any $j\in\mathcal{M}_*$, $\ell_j(0)$ will
%diverging to infinity in probability as $n\rightarrow\infty$. For
%finite sample situation, $\ell_j(0)$ will take large values for any
%$j\in\mathcal{M}_*$ with high probability. %Thus we propose a
%screening procedure based on the magnitude of $\ell_j(0)$.

%Aim to establish a general theoretical framework to cover both
%linear regression model and generalized model, $\ell_j(0)$ in this
%section means $\ell_j(0)$ in linear regression model and
%$\ell_j(0,0)$ in generalized linear model. Following two mild
%assumptions are need to guarantee our theoretical results.

An effective marginal screening procedure requires two conditions: (i)
if \mbox{$j\in\mathcal{M}_* $}, then $\rho_j$ takes nonnegligible
value; and (ii) if $j\notin\mathcal{M}_*$, then $\rho_j$ takes
negligible value. Actually, the first requirement is closely
related to recruiting the true signals that contribute to the response,
and the second one affects
the size of selected variable set that may contain false signals.
\citet
{FanLv2008JRSSB} shows that
under the identification condition
$\min_{j\in\mathcal{M}_*}|\rho_j|\geq f_n>0$ for some
function $f_n$, the first requirement is fulfilled. %In practice,
%$\rho_j$s are unknown and estimators $\hat{\beta}_j^M$s can be
%obtained by using (\ref{eqlrmg}). Subsequently, a screening
%procedure can be conducted by assessing the magnitudes of
%$\hat{\beta}_j^M$s.
%Generally speaking,
%$\hat{\beta}_j^M-\beta_j^M=O_p(n^{-{1}/{2}})$ under some regular
%conditions, which implies that $f_n$ should not weaken to zero
%faster than the rate $n^{-{1}/{2}}$. Hence,
A common assumption
for $f_n$ is $f_n=O(n^{-\kappa})$ for some
$\kappa\in(0,\frac{1}{2})$.

Our next theoretical analysis imposes the following two assumptions:

A.1: The random variable $Y$ has bounded
variance and there exists a positive constant $c_1$ such that
\[
\min_{j\in\mathcal{M}_*}\bigl|{\mathbb E}(X_j Y)\bigr|=\min
_{j\in\mathcal{M}_*}\bigl|{\operatorname{cov}(Y,X_{j})}\bigr|\geq
c_1n^{-\kappa}
\]
for some $\kappa\in[0,\frac{1}{2})$.

A.2: There are positive constants $K_1$, $K_2$, $\gamma_1$ and
$\gamma_2$ such that
\begin{eqnarray*}
\mathbb{P}\bigl\{|X_{j}|\geq u\bigr\}&\leq& K_1\exp
\bigl(-K_2u^{\gamma_1}\bigr) \qquad\mbox{for each } j=1,\ldots,p
\mbox{ and any } u>0,
\\
\mathbb{P}\bigl\{|Y|\geq u\bigr\}&\leq& K_1\exp
\bigl(-K_2u^{\gamma_2}\bigr) \qquad\mbox{for any } u>0.
\end{eqnarray*}

Assumption A.1 can be viewed as a requirement for the minimal signal
strength, and we call it the identification condition for
$j\in\mathcal{M}_*$. For linear models, the assumption
A.1 is same as condition 3 in \citet{FanLv2008JRSSB} that is
commonly assumed in sure independence feature screening.
For
generalized linear models, \citet{FanSong2010AOS} imposes the
identification condition as
${\min_{j\in\mathcal{M}_*}}|{\operatorname{cov}}(b'(\bX^{\mathrm{T}}\bbeta
),X_j)|\geq
c_1n^{-\kappa}$. %Actually,
By noticing that
$\operatorname{cov}(b'(\bX^{\mathrm{T}}\bbeta),X_j)=\mathbb{E}(X_jY)$, their
identification condition for $j\in\mathcal{M}_*$ is also same as
A.1. Since we impose no distributional assumptions, A.2 is assumed
to ensure the large deviation results that are used to get the
exponential convergence rate. The first part of A.2 is same as
the first part of condition D in \citet{FanSong2010AOS}. For linear
regression model, the second
part of condition D in \citet{FanSong2010AOS}
is equivalent to that $\bX^{\mathrm{T}}\bbeta$ satisfies the
Cram\'{e}r condition such that there exists a positive constant $H$
such that
$\mathbb{E}\{\exp(t\bX^{\mathrm{T}}\bbeta)\}<\infty$ for any
$|t|<H$. If the error
$\varepsilon$ is independent of covariates and satisfies the Cram\'{e}r
condition,
then we can obtain that the variable $Y$ also
satisfies the Cram\'{e}r condition. From Lemma 2.2 in
\citet{Petrov1995}, a random variable $W$ satisfies Cram\'{e}r
condition is
equivalent to that there are positive constants $b_1$ and $b_2$
such that $\mathbb{P}\{|W|\geq u\}\leq b_1\exp(-b_2u)$ for any
$u>0$. Therefore, our assumption here is actually weaker than that in
\citet{FanSong2010AOS}.
On the other hand, A.2 is also a general technical assumption in the literature
of large derivations. For example, $\gamma_1=2$ if $X_j$'s follow
normal distribution
or sub-Gaussian distribution, and $\gamma_1=\infty$ if $X_j$'s have
compact support. %}\footnote{discuss a few examples like the Gaussian.
%%in Proposition 2.

We now establish the following general result for the distribution
of empirical likelihood ratio which is the foundation for our
future theoretical results.
%
%th1 #&#
%
\begin{tm}\label{tm1}
Suppose\vspace*{1pt} that $U_1,\ldots,U_n$ are independent and identically
distributed random variables. Assume that there exist three positive
constants $\widetilde{K}_1$, $\widetilde{K}_2$ and $\gamma$ such
that $\mathbb{P}\{|U_i|>u\}\leq
\widetilde{K}_1\exp(-\widetilde{K}_2u^\gamma)$ for all $u>0$.
Define\vspace*{1pt}
$\mu_0=\mathbb{E}(U_i)$, $\delta=\max\{\frac{2}{\gamma}-1,0\}$,
$H=2^{1+\delta}$ and
$\bar{\Delta}=\frac{n^{{1/2}}\sigma}{2K}$, where
$\sigma^2=\mathbb{E}\{(U_i-\mu_0)^2\}$ and $K>\sigma$ is a
sufficiently\vspace*{1pt} large positive constant depending only on
$\widetilde{K}_1$, $\widetilde{K}_2$, $\gamma$ and $\mu_0$, then for
$L\rightarrow\infty$, there exists a positive constant $C$ only
depending on $\widetilde{K}_1$, $\widetilde{K}_2$ and $\gamma$ such
that
\[
\mathbb{P} \biggl\{\ell(\mu)<\frac{ n(\mu-\mu_0)^2}{L^2} \biggr\} \leq
\cases{\displaystyle \exp\biggl
\{-\frac{n(\mu-\mu_0)^2}{4H\sigma^2} \biggr\} +\exp\bigl(-CL^\gamma\bigr),
\vspace*{2pt}\cr
\qquad\mbox{if
}n^{{1/2}}|\mu-\mu_0|\leq\sigma\bigl(H^{1+\delta}\bar{
\Delta}\bigr)^{{1}/({1+2\delta})};
\vspace*{2pt}\cr
\displaystyle \exp\biggl\{-\frac{1}{4} \biggl(
\frac{n|\mu-\mu_0|}{2K} \biggr)^{{1}/({1+\delta})} \biggr\}+\exp\bigl
(-CL^\gamma
\bigr),
\vspace*{2pt}\cr
\qquad\mbox{if }n^{{1}/{2}}|\mu-\mu_0|> \sigma
\bigl(H^{1+\delta}\bar{\Delta}\bigr)^{{1}/({1+2\delta})};}
\]
where $\ell(\mu)$ is defined in Proposition \ref{pn1}.
\end{tm}

The proof of Theorem \ref{tm1} is given in \citet{Changetal2013},
where the main idea is applying large deviation theory
[\citet{SaulisStatuleviius1991,Petrov1995}].

Theorem \ref{tm1} reveals the magnitude of the empirical likelihood ratio
statistic evaluated at arbitrary values. When $\mu-\mu_0$ does not
diminish to 0, Theorem~\ref{tm1} implies that the empirical likelihood ratio
statistic diverges % takes diverging value
with large probability where the
diverging rate synthetically depends on the sample size $n$, some
diverging $L$ and the deviation of $\mu$ from the truth. Here $L$ is
a general technical device whose diverging rate is arbitrary. As a
direct result of Theorem \ref{tm1}, we have the following proposition for
$\ell_j(0)$.
%the tail probability. In fact, we only need to require $X_{ij}Y_i$
%satisfies the Cram\'{e}r condition. Note that
%the last inequality is based on the fact
%$|\mathbb{E}X_{ij}Y_i|=|\operatorname{cov}(X_{ij},Y_i)|\leq\sigma_Y$
%where $\sigma_Y^2=\operatorname{var}(Y)$. If $u\leq\sigma_Y$, we can
%choose $K_1$ and $K_2$ such that $K_1\geq\exp(K_2\sigma_Y)$. If
%$u>\sigma_Y$, A.2 follows the condition that
%$X_{ij}Y_i$ satisfies the Cram\'{e}r condition. On the other hand,
%for any $v>0$,
%If the tail probabilities of $X_{ij}$ and $Y_i$ are exponential
%decreasing with the parameter $\alpha=2$ in
%then A.2 holds.
%

%
%pr2 #&#
%
\begin{pn}\label{prop2}
Under assumptions \textup{A.1} and \textup{A.2}, there exists a positive constant
$C_1$ depending only on $K_1$, $K_2$, $\gamma_1$ and $\gamma_2$
appeared in assumption \textup{A.2} such that for any $j\in\mathcal{M}_*$ and
$L\rightarrow\infty$,
\[
\mathbb{P} \biggl\{\ell_j(0)<\frac{c_1^2
n^{1-2\kappa}}{L^2} \biggr\}\leq\cases{
\exp\bigl(-C_1n^{1-2\kappa}\bigr)+\exp\bigl(-C_1L^{\gamma}
\bigr),
\vspace*{2pt}\cr
\qquad\mbox{if }(1-2\kappa) (1+2\delta)<1;
\vspace*{2pt}\cr
\exp\bigl(-C_1n^{({1-\kappa})/({1+\delta})}
\bigr)+\exp\bigl(-C_1L^{\gamma}\bigr),
\vspace*{2pt}\cr
\qquad\mbox{if }(1-2\kappa)
(1+2\delta)\geq1;}
\]
where $\gamma=\frac{\gamma_1\gamma_2}{\gamma_1+\gamma_2}$ and
$\delta=\max\{\frac{2}{\gamma}-1,0\}$.
\end{pn}

Proposition \ref{prop2} is a uniform result for all features
contributing in
the true model. Specifically, with large probability and uniformly
for all $j\in\mathcal{M}_*$, the diverging rate of $\ell_j(0)$ is
not slower than $n^{1-2\kappa}L^{-2}$. From Proposition \ref{pn1}, if
$|\mathbb{E}(X_jY)|=O(n^{-w})$ for some $w\in(0,\frac{1}{2})$ and some
$j\in\mathcal{M}_*$, then $\ell_j(0)\xrightarrow{p}\infty$. This
can be viewed as a requirement such that the signal strength cannot
diminish to 0 at a too fast rate. Therefore,
$n^{{1}/{2}-\kappa}L^{-1}\rightarrow\infty$ as $n\rightarrow
\infty$ is
required for sure independence screening. By choosing
$L=n^{{1}/{2}-\kappa-\tau}$ for some $\tau\in(0,\frac
{1}{2}-\kappa)$, we obtain
the following corollary more specifically summarizing that the set
$\mathcal{M}_*$ can be distinguished by examining the marginal
empirical likelihood ratio $\ell_j(0)$ $(j=1,\ldots,p)$.
%
%co1 #&#
%
\begin{cy}
Under assumptions \textup{A.1} and \textup{A.2}, there exists a positive constant
$C_1$ depending only on $K_1$, $K_2$, $\gamma_1$ and $\gamma_2$
appeared in assumption~\textup{A.2} such that, for any
$\tau\in(0,\frac{1}{2}-\kappa)$,
\[
\max_{j\in\mathcal{M}_*}\mathbb{P}\bigl\{\ell_j(0)<c_1^2
n^{2\tau}\bigr\} \leq\cases{\exp\bigl\{-C_1n^{(1-2\kappa)\wedge
({(1-2\kappa-2\tau
)\gamma}/{2})}
\bigr\},
\vspace*{2pt}\cr
\qquad\mbox{if }(1-2\kappa) (1+2\delta)<1;
\vspace*{2pt}\cr
\exp\bigl
\{-C_1n^{(({1-\kappa})/({1+\delta}))\wedge({(1-2\kappa
-2\tau)\gamma}/{2})}\bigr\},
\vspace*{2pt}\cr
\qquad\mbox{if }(1-2\kappa) (1+2\delta)
\geq1;}
\]
where $\gamma=\frac{\gamma_1\gamma_2}{\gamma_1+\gamma_2}$ and
$\delta=\max\{\frac{2}{\gamma}-1,0\}$.
\end{cy}

Summarizing above results, we formally establish the screening
properties of the marginal empirical likelihood approach.
%
%th2 #&#
%
\begin{tm}\label{tm2}
Under assumptions \textup{A.1} and \textup{A.2}, there exists a positive constant
$C_1$ depending only on $K_1$, $K_2$, $\gamma_1$ and $\gamma_2$
appeared in assumption \textup{A.2} such that, for any
$\tau\in(0,\frac{1}{2}-\kappa)$ and $\gamma_n=c_1^2 n^{2\tau}$,
\begin{eqnarray*}
\mathbb{P}\{\mathcal{M}_*\subset\widehat{\mathcal{M}}_{\gamma_n}\} \geq
\cases{1-s\exp\bigl\{-C_1n^{(1-2\kappa)\wedge({(1-2\kappa-2\tau
)\gamma}/{2})}\bigr\},
\vspace*{2pt}\cr
\qquad\mbox{if }(1-2
\kappa) (1+2\delta)<1;
\vspace*{2pt}\cr
1-s\exp\bigl\{-C_1n^{(({1-\kappa})/({1+\delta}))\wedge({(1-2\kappa
-2\tau)\gamma}/{2})}\bigr\},
\vspace*{2pt}\cr
\qquad\mbox{if }(1-2\kappa) (1+2\delta)\geq1;}
\end{eqnarray*}
where $\gamma=\frac{\gamma_1\gamma_2}{\gamma_1+\gamma_2}$ and
$\delta=\max\{\frac{2}{\gamma}-1,0\}$.
\end{tm}
Theorem \ref{tm2} implies the sure screening property for our
procedure with nonpolynomial dimensionality:
\[
\log p=\cases{ o\bigl(n^{(1-2\kappa)\wedge({(1-2\kappa-2\tau)\gamma
}/{2})}\bigr), &\quad if $(1-2\kappa) (1+2\delta)<1$;
\vspace*{2pt}\cr
o\bigl(n^{(({1-\kappa})/({1+\delta}))\wedge({(1-2\kappa-2\tau
)\gamma}/{2})}\bigr), &\quad if $(1-2\kappa) (1+2\delta)\geq1$.}
\]
When the covariates and error are normal, $\gamma_1=2$ and
$\gamma_2=2$. Then $\gamma=1$, $\delta=1$ and $\log
p=o(n^{{1}/{2}-\kappa})$ which is weaker than that in \citet
{FanLv2008JRSSB} where $\log p = o(n^{1-2\kappa})$ is allowed.
This can be viewed as a price paid for allowing nonnormal covariate
and more general error distribution.
%However, our assumptions are significantly weaker than theirs.
Furthermore, we compare our result and that in \citet{FanSong2010AOS}.
The Lemma 1 in \citet{FanSong2010AOS} means that $\gamma_2=1$. The
corresponding parameters under their this setting are
$\gamma=\frac{\gamma_1}{\gamma_1+1}$ and
$\delta=\frac{\gamma_1+2}{\gamma_1}$, respectively. Then, we can
handle the nonpolynomial dimensionality
\[
\log p=o\bigl(n^{{(1-2\kappa)\gamma_1}/({2\gamma_1+2})}\bigr)
\]
in this setting, which is actually a stronger result than that in
\citet
{FanSong2010AOS} where $\log p=o(n^{{(1-2\kappa)\gamma_1}/{A}})$ and
$A=\operatorname{max}\{\gamma_1+4,3\gamma_1+2\}$.\vspace*{1pt}

Now we investigate how large the set
$\widehat{\mathcal{M}}_{\gamma_n}$ is. This question is closely
related to the asymptotic property of $\ell_j(0)$ for
$j\notin\mathcal{M}_*$. %Let
Essentially, we need to know the magnitudes of $\ell_j(0)$ for
$j\notin\mathcal{M}_*$. We first consider the simple case
$\rho_j=0$
for any $j\notin\mathcal{M}_*$ and have the following result. %

%
%pr3 #&#
%
\begin{pn}\label{prop3}
Under assumptions \textup{A.1} and \textup{A.2}, if $\rho_j=0$,
there is a positive constant $C_2$ depending\vspace*{1pt} only on $K_1$, $K_2$,
$\gamma_1$ and $\gamma_2$ appeared in assumption~\textup{A.2} such that, for
any $\tau\in(0,\frac{1}{2}-\kappa)$,
\[
\mathbb{P}\bigl\{\ell_j(0)\geq c_1^2n^{2\tau}
\bigr\}\leq\cases{ \exp\bigl(-C_2n^{2\tau}\bigr), &\quad if $
\gamma<4$ and $\displaystyle \tau\leq\frac
{\gamma}{12}$;
\vspace*{2pt}\cr
\exp\bigl(-C_2n^{{\gamma}/{6}}
\bigr), &\quad if $\gamma<4$ and $\displaystyle \tau>\frac{\gamma}{12}$;
\vspace*{2pt}\cr
\exp
\bigl(-C_2n^{2\tau}\bigr), &\quad if $\gamma\geq4$ and $\displaystyle \tau
\leq\frac{\gamma}{2\gamma
+4}$;
\vspace*{2pt}\cr
\exp\bigl(-C_2n^{{\gamma}/({\gamma+2})}\bigr), &\quad if
$\gamma\geq4$ and $\displaystyle \tau>\frac{\gamma}{2\gamma+4}$;}
\]
where $\gamma=\frac{\gamma_1\gamma_2}{\gamma_1+\gamma_2}$.
\end{pn}

The assumption $\rho_j=0$ for any $j\notin\mathcal{M}_*$ can be
guaranteed by the partial orthogonality condition, that is,
$\{X_j\dvtx j\notin\mathcal{M}_*\}$ is independent of
$\{X_j\dvtx j\in\mathcal{M}_*\}$. The orthogonality condition is
essentially the assumption made in \citet{HuangHortowitzMa2008AOS}
who showed the model selection consistency in the case with the
ordinary linear model and bridge regression. This proposition gives
the property of $\ell_j(0)$ for any $j\notin\mathcal{M}_*$ which can
be used to establish the theoretical result for the size of
$\widehat{\mathcal{M}}_{\gamma_n}$ where $\gamma_n=c_1^2n^{2\tau}$.
Note that\looseness=-1
\begin{eqnarray*}
|\widehat{\mathcal{M}}_{\gamma_n}|&=&\sum_{j\in\mathcal{M}_*}I
\bigl\{ \ell_j(0)\geq c_1^2n^{2\tau}
\bigr\}+\sum_{j\notin\mathcal{M}_*}I\bigl\{\ell_j(0)\geq
c_1^2n^{2\tau}\bigr\}
\\[-2pt]
&\leq& s+\sum_{j\notin\mathcal{M}_*}I\bigl\{\ell_j(0)\geq
c_1^2n^{2\tau}\bigr\},
\end{eqnarray*}\looseness=0
then
\[
\mathbb{P}\bigl\{|\widehat{\mathcal{M}}_{\gamma_n}|>s\bigr\} \leq\sum
_{j\notin\mathcal{M}_*}\mathbb{P}\bigl\{\ell_j(0)\geq
c_1^2n^{2\tau}\bigr\}.
\]
By Proposition \ref{prop3}, we obtain the following theorem.

%
%th3 #&#
%
\begin{tm}\label{tm3}
Under assumptions \textup{A.1} and \textup{A.2}, if $\rho_j=0$ for any $j\notin\mathcal
{M}_*$, %for $\tau\in(0,1/2-\kappa)$
then there exists a positive constant $C_2$ depending only on $K_1$,
$K_2$, $\gamma_1$ and $\gamma_2$ appeared in assumption \textup{A.2} such
that, for any $\tau\in(0,\frac{1}{2}-\kappa)$ and
$\gamma_n=c_1^2n^{2\tau}$,
\[
\mathbb{P}\bigl\{|\widehat{\mathcal{M}}_{\gamma_n}|>s\bigr\} \leq\cases{ p\exp
\bigl(-C_2n^{2\tau}\bigr),&\quad if $\gamma<4$ and $\tau\leq
\displaystyle \frac
{\gamma}{12}$;
\vspace*{1pt}\cr
p\exp\bigl(-C_2n^{{\gamma}/{6}}\bigr), &\quad if $
\gamma<4$ and $\displaystyle \tau>\frac{\gamma}{12}$;
\vspace*{1pt}\cr
p\exp\bigl(-C_2n^{2\tau}
\bigr), &\quad if $\gamma\geq4$ and $\displaystyle \tau\leq\frac{\gamma}{2\gamma
+4}$;
\vspace*{1pt}\cr
p\exp
\bigl(-C_2n^{{\gamma}/({\gamma+2})}\bigr), &\quad if $ \gamma\geq4$
and $\displaystyle\tau>\frac{\gamma}{2\gamma+4}$;}
\]
where $\gamma=\frac{\gamma_1\gamma_2}{\gamma_1+\gamma_2}$.
\end{tm}
%
%The specification of the positive constant $\xi$ are given in the
%Appendix.

%This assumption is also
%satisfied in the setting of Fan and Lv (2008) for linear regression
%model. In their setting, $\gamma_1=\gamma_2=2$ which imply
%$\gamma=1$.
From Theorem \ref{tm3}, we have $
\mathbb{P}\{|\widehat{\mathcal{M}}_{\gamma_n}|>s\}\leq
p\exp\{-C_2n^{(2\tau)\wedge({\gamma}/{6})\wedge({\gamma
}/({\gamma+2}))}\}$
which means that the event $\{|\widehat{\mathcal{M}}_{\gamma_n}|\leq
s\}$ occurs with probability approaching to $1$ if $\log
p=o(n^{(2\tau)\wedge({\gamma}/{6})\wedge({\gamma
}/({\gamma+2}))})$.
On the other hand, following Theorem \ref{tm2}, we have
%Theorem 2, we have that
$\mathbb{P}\{\mathcal{M}_*\subset\widehat{\mathcal{M}}_{\gamma
_n}\}\rightarrow1$
provided $\log
p=o(n^{({(1-2\kappa-2\tau)\gamma}/{2})\wedge(1-2\kappa)})$.
Combining these two results together, we can obtain that
\[
\mathbb{P}\{\widehat{\mathcal{M}}_{\gamma_n}=\mathcal{M}_*\} \rightarrow1
\qquad\mbox{if } \log p=o\bigl(n^{({\gamma}/{6})\wedge({(1-2\kappa
)\gamma}/({\gamma
+2}))}\bigr)
\]
and
\[
\tau=
\frac{(1-2\kappa)\gamma}{2\gamma+4}.
\]
This property shows the selection consistency of our procedure. In a
more general case without partial orthogonality\vadjust{\goodbreak} condition, we could
consider the size of the set
$\widehat{\mathcal{M}}_{\gamma_n}$ under the setting
\[
\max_{j\notin\mathcal{M}_*}|\rho_j|=o\bigl(n^{-\kappa}\bigr),
\]
which is an assumption imposed in \citet{FanSong2010AOS}. %
%in Theorem 3, $\xi_1$ and $\xi_2$ in Theorem 4 in the Appendix. }

%
%pr4 #&#
%
\begin{pn}\label{prop4}
Under assumptions \textup{A.1} and \textup{A.2}, if
${\max_{j\notin\mathcal{M}_*}}|\rho_j|=\break O(n^{-\eta})$ where
$\eta>\kappa$ and
$\min_{j\notin\mathcal{M}_*}\mathbb{E}(X_j^2Y^2)\geq c_2$ for some
$c_2>0$, there exists a positive constant $C_3$ depending only on
$K_1$, $K_2$, $\gamma_1$ and $\gamma_2$ appeared in assumption
\textup{A.2}
and $c_2$ such that, for any $j\notin\mathcal{M}_*$ and
$\tau\in(\frac{1}{2}-\eta,\frac{1}{2}-\kappa)$,
{\fontsize{10pt}{12pt}\selectfont{
\begin{eqnarray*}
&&\mathbb{P}\bigl\{\ell_j(0)\geq c_1^2n^{2\tau}
\bigr\}
\\[-2pt]
&&\qquad\leq\cases{ \exp\bigl(-C_3n^{2\tau}\bigr)+\exp
\bigl(-C_3n^{{\gamma}/{6}}\bigr),&\quad if $\gamma<2$ and $\eta>
\displaystyle \frac{1}{4}$;
\vspace*{2pt}\cr
\exp\bigl(-C_3n^{\gamma\eta}\bigr)+\exp
\bigl(-C_3n^{{\gamma}/{6}}\bigr), &\quad if $\gamma<2$ and $\displaystyle \eta
\leq\frac{1}{4}$;
\vspace*{2pt}\cr
\exp\bigl(-C_3n^{\gamma\eta}\bigr)+
\exp\bigl(-C_3n^{{\gamma}/{6}}\bigr), &\quad if $\gamma\geq2$ and
$\displaystyle \eta\leq\frac{1}{\gamma+2}$;
\vspace*{2pt}\cr
\exp\bigl(-C_3n^{{\gamma}/({\gamma+2})}\bigr)+
\exp\bigl(-C_3n^{2\tau}\bigr), &\quad if $\gamma\geq4$ and
$\displaystyle \eta>\frac{1}{\gamma+2}$;
\vspace*{2pt}\cr
\exp\bigl(-C_3n^{{\gamma}/{6}}\bigr)+
\exp\bigl(-C_3n^{2\tau}\bigr), &\quad if $2\leq\gamma<4$
and $\displaystyle \eta>\frac{1}{\gamma+2}$;}
\end{eqnarray*}}}

\noindent where $\gamma=\frac{\gamma_1\gamma_2}{\gamma_1+\gamma_2}$.
\end{pn}

If $\rho_j=0$ for any $j\notin\mathcal{M}_*$, then $\eta=\infty$.
Hence, this proposition reduces to Proposition \ref{prop3}. Following
the same argument between Proposition \ref{prop3} and
Theorem~\ref{tm3}, we can obtain the following theorem related to the
size of $\widehat{\mathcal{M}}_{\gamma_n}$.

%
%th4 #&#
%
\begin{tm} \label{tm4}
Under assumptions \textup{A.1} and \textup{A.2}, if
$\max_{j\notin\mathcal{M}_*}|\rho_j|=O(n^{-\eta})$ where
$\eta>\kappa$ and
$\min_{j\notin\mathcal{M}_*}\mathbb{E}(X_j^2Y^2)\geq c_2$ for some
$c_2>0$, then there exists a positive constant
$C_3$ only
depending on $K_1$, $K_2$, $\gamma_1$ and $\gamma_2$ appeared in
assumption~\textup{A.2} and $c_2$ such that, for any
$\tau\in(\frac{1}{2}-\eta,\frac{1}{2}-\kappa)$ and
$\gamma_n=c_1^2n^{2\tau}$,
%Pick $\tau\in(1/2-\eta,1/2-\kappa)$,
%
{\fontsize{10pt}{12pt}\selectfont{
\begin{eqnarray*}
\hspace*{-4pt}&&\mathbb{P}\bigl\{|\widehat{\mathcal{M}}_{\gamma_n}|>s\bigr\}
\\[-2pt]
\hspace*{-4pt}&&\qquad\leq\cases{ p\exp\bigl(-C_3n^{2\tau}\bigr)+p\exp
\bigl(-C_3n^{{\gamma}/{6}}\bigr), &\quad if $\gamma<2$ and $\eta>
\displaystyle \frac{1}{4}$;
\vspace*{2pt}\cr
p\exp\bigl(-C_3n^{\gamma\eta}\bigr)+p\exp
\bigl(-C_3n^{{\gamma}/{6}}\bigr), &\quad if $\gamma<2$ and $\displaystyle \eta
\leq\frac{1}{4}$;
\vspace*{2pt}\cr
p\exp\bigl(-C_3n^{\gamma\eta}\bigr)+p
\exp\bigl(-C_3n^{{\gamma}/{6}}\bigr), &\quad if $\gamma\geq2$ and
$\displaystyle \eta\leq\frac{1}{\gamma+2}$;
\vspace*{2pt}\cr
p\exp\bigl(-C_3n^{{\gamma}/({\gamma+2})}
\bigr)+p\exp\bigl(-C_3n^{2\tau}\bigr), &\quad if $\gamma\geq4$
and $\displaystyle \eta>\frac{1}{\gamma+2}$;
\vspace*{2pt}\cr
p\exp\bigl(-C_3n^{{\gamma}/{6}}
\bigr)+p\exp\bigl(-C_3n^{2\tau}\bigr), &\quad if $2\leq\gamma<4$
and $\displaystyle \eta>\frac{1}{\gamma+2}$;}
\end{eqnarray*}}}

\noindent where
$\gamma=\frac{\gamma_1\gamma_2}{\gamma_1+\gamma_2}$.\vadjust{\goodbreak}
\end{tm}

In summary, our results show that the marginal empirical likelihood
approach has a very good control of the size of the recruited
variables. With large probability, the set of the recruited
variables is not larger than the true contributing explanatory
variables. As shown later in our simulation results, the marginal
empirical likelihood approach perform very well in terms of the set
of false selected variables by the marginal empirical likelihood
approach.

%s4 #&#
\section{Extensions} \label{s4}

%s4.1 #&#
\subsection{A broad framework}

The marginal empirical likelihood %and its iterative version
can be applied in a general framework besides the linear
models and generalized linear models. Based on general estimating
equations approach [\citet{Hansen1982ecnoca,QinLawless1994AS}], we
can also apply the screening procedure based on the marginal
empirical likelihood. We will demonstrate that the marginal
empirical likelihood approach provides an effective device to
combine information that can be used to enhance the performance of a
screening procedure.

%{\bf Here, I transfer $\theta$ to $\beta$ in order to have the same
%parameter form as the example.}

Let $\bZ_i\in{\mathbb R}^d$ $(i=1,\ldots,n)$ be generic
observations, $\bbeta=(\beta_1,\ldots,\beta_p)^{\mathrm{T}}\in
{\mathbb R}^p$
be parameter of interest and
$\bg(\bZ;\bbeta)=(g_1(\bZ;\bbeta),\ldots,g_r(\bZ;\bbeta
))^{\mathrm{T}}$
be the
$r$-dimensional estimating function such that ${\mathbb
E}\{\bg(\bZ;\bbeta)\}=\bzero$. Let $\mathcal{M}_*=\{1\leq j\leq
p\dvtx\beta_{j}\neq0\}$ be the true model with size $|\mathcal{M}_*|=s$.
We are interested in how to construct a sure feature screening
procedure to recover $\mathcal{M}_*$ in the general estimating
equation
setting. %Similar to the previous arguments of marginal empirical
%likelihood approach proposed for linear models and
%generalized linear models, we consider the
To motivate the marginal empirical likelihood approach, let us consider
the estimating function evaluated at
\[
\bbeta^{(j)}=(\underbrace
{0,\ldots,0}_{j-1},\beta,\underbrace{0,\ldots,0}_{p-j})^{\mathrm{T}}\qquad
(j=1,\ldots, p).
\]
In practice, many components in $\bg(\bZ;\bbeta^{(j)})$ do not
involve the unknown parameter; see, for example, the estimating
function constructed from the least-squares method and our example
given later.
%marginal estimating function evaluated for the $j$th component of $
%$
%$
%Ignoring some
%equations in $\bg(\bZ;(0,\ldots,0,\beta,0,\ldots,0))$ that are exact
%equal to zero, we rewrite $\bg^{(j)}(\bZ;\beta)$ as
Therefore, we denote by
\[
\bg^{(j)}(\bZ;\beta) =\bigl(g_1^{(j)}(\bZ;\beta),\ldots,g^{(j)}_{r_j}(\bZ;\beta)\bigr)^{\mathrm{T}}
\]
an $r_j$ ($r_j\ge1$)-dimensional estimating function collecting the
components in $\bg(\bZ;\bbeta^{(j)})$ that depend on the unknown
parameter. Usually $r_j>1$ is small and not all components of $\bZ$ are
involved in $\bg^{(j)}(\bZ;\beta)$.
A remarkable advantage of this broad framework is that it provides a
device for feature screening using
more flexibly constructed conditions so that additional
data information can be more effectively incorporated.

Correspondingly, we define the marginal empirical
likelihood for $\beta$ as
%
%e4.1 #&#
%
\begin{equation}\label{eqel2}
\mathrm{EL}_j(\beta) =\sup\Biggl\{ \prod
_{i=1}^n w_i\dvtx w_i\ge0,
\sum_{i=1}^n w_i=1, \sum
_{i=1}^n w_i
\bg^{(j)}(\bZ_i;\beta)=\bzero\Biggr\}.
\end{equation}
Then screening can be done based on the ranking of $\mathrm{EL}_j(0)$ or
equivalently using the corresponding marginal empirical likelihood
ratio evaluated at 0---that is, $\ell_j(0)$.
The steps of the procedure are the same as those described earlier.
%and the iterative screening can also be extended.
A concrete example of this scenario is given as follows.

\begin{example*}[{(Quadratic inference function (QIF) approach
[\citet{Quetal2000Bioka}])}]
Longitudinal data arise commonly in biomedical research with
repeated measurements from the same subject or within the same
cluster. Let $Y_{it}$ and $\bX_{it}$ $(i=1,\ldots, n,
t=1,\ldots,m_i)$ be the response and covariates of the $i$th subject
measured at time $t$. Let $\mathbb{E}(Y_{it})=\mu(\bX^{\mathrm{T}}
_{it}\bbeta)=\mu_{it}$ where $\bbeta\in{\mathbb R}^p$ is the
parameter of interest. Incorporating the dependence among the
repeated measurements is essential for efficient inference. \citet
{LiangZeger1986Bioka} proposed to estimate $\bbeta$ by solving
$0=\sum_{i=1}^n \dot{\bmu}_i ^{\mathrm{T}}\bW_i^{-1}(\bY_i-\bmu
_i)$. %
Here for the $i$th subject, $\bY_i=(Y_{i1},\ldots, Y_{i
m_i})^{\mathrm{T}}$,
$\bmu_i=(\mu_{i1},\ldots, \mu_{im_i})^{\mathrm{T}}$, $\dot{\bmu
}_i=\frac
{\partial\bmus_i}{\partial\bbetaa} $ and $\bW_i=\bv_i^{
{1}/{2}} \bR\bv_i^{{1}/{2}}$,
where $\bv_i$ is a diagonal matrix of the conditional variances of
subject $i$ and $\bR$ is
a working correlation matrix that may depend on some unknown parameter.
This approach uses estimating function $\bg(\bZ_i;\bbeta)=\dot{\bmu
}_i ^{\mathrm{T}}\bW_i^{-1}(\bY_i-\bmu_i)$,
where $\bZ_i=(\bZ_{i1}^{\mathrm{T}},\ldots,\bZ_{im_i}^{\mathrm
{T}})^{\mathrm{T}}$,
$\bZ_{it}=(Y_{it},\bX_{it}^{\mathrm{T}})^{\mathrm{T}}$ and $r=p$. %
%dispersion parameter by the method of moments.
More recently, \citet{Quetal2000Bioka} proposed to model $\bR
^{-1}$ by
$\sum_{i=1}^m a_i \bM_i$, where $\bM_1,\ldots, \bM_m$ are known
matrices and $a_1,\ldots, a_m$ are unknown constants. Then $\bbeta$
can be estimated by the quadratic inference functions approach
[\citet{Quetal2000Bioka}] that uses
%
%e4.2 #&#
%
\begin{equation}\label{eqgee2}
\bg(\bZ_i;\bbeta)=\pmatrix{\dot{\bmu}_i
^{\mathrm{T}}\bv_i^{-
{1}/{2}}\bM_1
\bv_i^{-{1}/{2}}(\bY_i-\bmu_i)
\cr
\vdots
\cr
\dot{\bmu}_i ^{\mathrm{T}}\bv_i^{-{1}/{2}}
\bM_m \bv_i^{-{1}/{2}}(\bY_i-
\bmu_i)}\qquad(i=1,\ldots,n).
\end{equation}
This falls into our framework with $r>p$ when $m>1$, and with
$r=p$ if $m=1$. When applying the marginal approach, we note that
$\bg^{(j)}(\bZ;\beta)$ is an $m$-dimensional estimating function.
The marginal screening by empirical likelihood can be conveniently
applied to this scenario, and we note that the existing
independence screening methods cannot be directly applied when
$m>1$.

In a concurrent and independent work, \citet{ZhaoLi2012}
considered
feature screening using estimating\vspace*{1pt} functions when $r=p$. By using our
notations, their approach are based on $\bg^{(j)}(\bZ;0)$---the
marginal estimating function evaluated at~0. Their screening
procedure are based on ranking the %quantities containing the
%respective component in each
absolute value of $\bg^{(j)}(\bZ;0)$ %with the largest
%absolute value
for $j=1,\ldots, p$. Our approach is different as seen
from the above marginal empirical likelihood construction. In
addition, analogous to that in linear models and
generalized linear models, the marginal empirical likelihood
constructed from using the marginal estimating function is also
capable of incorporating the level of uncertainties associated with
finite sample estimating functions.

We now characterize the properties of the screening procedure in the
framework of models specified by estimating equations. For any
vector ${\mathbf a}=(a_1,\ldots,a_q)^{\mathrm{T}}\in\mathbb{R}^q$, we
use $\|{\mathbf
a}\|_\infty=\max_{i=1,\ldots,q}|a_q|$ and $\|{\mathbf
a}\|_2=(\sum_{i=1}^qa_i^2)^{{1}/{2}}$ to denote its $L_\infty$
and $L_2$ norms, respectively. Aiming to establish the theoretical
results, we need the following two assumptions.

A.3: There exists a positive constant $c_3$ such that
\[
\min_{j\in\mathcal{M}_*}\bigl\|\mathbb{E}\bigl\{\bg^{(j)}(\bZ;0)\bigr\}
\bigr\|_\infty\geq c_3n^{-\kappa}
\]
for some $\kappa\in[0,\frac{1}{2})$.

A.4: There are positive constants $K_3$, $K_4$ and $\gamma_3$ such
that
\[
\mathbb{P}\bigl\{\bigl\|\bg^{(j)}(\bZ;0)\bigr\|_2\geq u\bigr\}\leq
K_3\exp\bigl(-K_4u^{\gamma_3}\bigr)
\]
for each $j=1,\ldots,p$ and any $u>0$.

Assumption A.3 is a general identification condition for
the set $\mathcal{M}_*$ %in the generalized estimating equations
%setting.
when considering the broad framework of models specified by general
estimating equations.
It means that the weakest signals reflected by
%there exists a component of
$\|\mathbb{E}\{\bg^{(j)}(\bZ;0)\}\|_\infty$ $(j\in\mathcal{M}_*)$
cannot vanish at a rate faster % such that its strength is stronger
than $n^{-{1}/{2}}$. Assumption A.3 is not stringent, and it
reduces to A.1 in special cases of linear models and
generalized linear models.
A similar assumption is also made in \citet{ZhaoLi2012}. Assumption
A.4, which is a counterpart of A.2 in general cases, is required for
%%stating the
%exponential inequality.
establishing exponential inequality when analyzing large deviations.
\citet{ZhaoLi2012} assumed
boundness of all components in $\bg^{(j)}(\bZ;0)$, which implies A.4.
\end{example*}

%
%th5 #&#
%
\begin{tm}\label{tm5}
Under assumptions \textup{A.3--A.4}, there exists a positive constant $C_4$
depending only on $K_3$, $K_4$ and $\gamma_3$ appeared in assumption
\textup{A.4} such that, for any $\tau\in(0,\frac{1}{2}-\kappa)$ and
$\gamma_n=c_3^2n^{2\tau}$,
\begin{eqnarray*}
\mathbb{P}\{\mathcal{M}_*\subset\widehat{\mathcal{M}}_{\gamma_n}\} \leq
\cases{ 1-s\exp\bigl\{-C_4n^{(1-2\kappa)\wedge({(1-2\kappa-2\tau
)\gamma
_3}/{2})}\bigr\},
\vspace*{2pt}\cr
\qquad\mbox{if }(1-2
\kappa) (1+2\delta)<1;
\vspace*{2pt}\cr
1-s\exp\bigl\{-C_4n^{(({1-\kappa})/({1+\delta}))\wedge({(1-2\kappa
-2\tau)\gamma_3}/{2})}\bigr\},
\vspace*{2pt}\cr
\qquad\mbox{if }(1-2\kappa) (1+2\delta)\geq1;}
\end{eqnarray*}
where $\delta=\max\{\frac{2}{\gamma_3}-1,0\}$.
\end{tm}

This theorem %can be viewed as the
is a natural
extension of Theorem \ref{tm2} in the
%generalized estimating equation setting.
broad framework for models specified by general estimating equations.
In %the screening procedure
special cases, we have considered for linear models and
generalized linear models, $\bg^{(j)}(\bZ;0)=X_jY$ $(j=1,\ldots,p)$,
and $\gamma_3$ in assumption A.4 is equal to
$\frac{\gamma_1\gamma_2}{\gamma_1+\gamma_2}$ where $\gamma_1$ and
$\gamma_2$ are %proposed for the distributions of $X_j$ and $Y$ in
specified in
A.2.

Let
\[
\bu_j=\mathbb{E}\bigl\{\bg^{(j)}(\bZ;0)\bigr\}
\qquad\mbox{for each }j=1,\ldots,p.
\]
We now consider the size of $\widehat{\mathcal{M}}_{\gamma_n}$ in the
setting
\[
{\max_{j\notin\mathcal{M}_*}}\|\bu_j\|_\infty=o
\bigl(n^{-\kappa}\bigr).
\]
This specification also reduces to those considered in special cases
of linear models and generalized linear models. The
counterpart of Theorem \ref{tm4} for establishing the selection consistency
is given as follows.

%and the corresponding theoretical result is proposed in the
%following theorem.
%
%th6 #&#
%
\begin{tm}\label{tm6}
Under assumptions \textup{A.3} and \textup{A.4}, if
${\max_{j\notin\mathcal{M}_*}}\|\bu_j\|_\infty=O(n^{-\eta})$ where
$\eta>\kappa$ and
$\min_{j\notin\mathcal{M}_*}\lambda_{\min}(\mathbb{E}\{\bg
^{(j)}(\bZ;0)\bg^{(j)}(\bZ;0)^{\mathrm{T}}\})\geq
c_4$ for some $c_4>0$, where $\lambda_{\min}(A)$ means the smallest
eigenvalue of $A$, then there exists a positive constant
$C_5$ depending only on $K_3$, $K_4$ and $\gamma_3$ appeared in
assumption \textup{A.4} and $c_4$ such
that, for any $\tau\in(\frac{1}{2}-\eta,\frac{1}{2}-\kappa)$ and
$\gamma_n=c_3n^{2\tau}$,
%Pick $\tau\in(1/2-\eta,1/2-\kappa)$,
%
\begin{eqnarray*}
\hspace*{-4pt}&&\mathbb{P}\bigl\{|\widehat{\mathcal{M}}_{\gamma_n}|>s\bigr\}
\\
\hspace*{-4pt}&&\qquad\leq\cases{ p\exp\bigl(-C_5n^{2\tau}\bigr)+p\exp
\bigl(-C_5n^{{\gamma_3}/{6}}\bigr), \qquad \mbox{if $\gamma_3<2$
and $\displaystyle \eta>\frac{1}{4}$};
\vspace*{2pt}\cr
p\exp\bigl(-C_5n^{\gamma_3\eta}
\bigr)+p\exp\bigl(-C_5n^{{\gamma_3}/{6}}\bigr), \hspace*{-4.01pt}\qquad \mbox{if $
\gamma_3<2$ and $\displaystyle \eta\leq\frac{1}{4}$};
\vspace*{2pt}\cr
p\exp
\bigl(-C_5n^{\gamma_3\eta}\bigr)+p\exp\bigl(-C_5n^{{\gamma_3}/{6}}
\bigr), \hspace*{-3.95pt}\qquad \mbox{if $\gamma_3\geq2$ and $\displaystyle \eta\leq\frac{1}{\gamma_3+2}$};
\vspace*{2pt}\cr
p\exp\bigl(-C_5n^{{\gamma_3}/({\gamma_3+2})}\bigr)+p\exp\bigl
(-C_5n^{2\tau}
\bigr), \cr
\hspace*{162pt}\qquad\mbox{if $\gamma_3\geq4$ and $\displaystyle \eta>\frac{1}{\gamma_3+2}$};
\vspace*{2pt}\cr
p\exp\bigl(-C_5n^{{\gamma_3}/{6}}\bigr)+p\exp\bigl(-C_5n^{2\tau}
\bigr), \qquad \mbox{if $2\leq\gamma_3<4$ and $\displaystyle \eta>\frac{1}{\gamma_3+2}$}.}
\end{eqnarray*}
\end{tm}

%This result is totally the same as the one proposed in Theorem 4.
Combining the Theorems \ref{tm5} and \ref{tm6}, we
can see that the screening procedure using the marginal empirical
likelihood ratio is valid in a broad framework for identifying the set
of the effective features.
%can find our screening procedure
%based on marginal empirical likelihood ratio is less restrictive on
%the distributional assumption, and is adapted to be applied in a
%broad range of problems.

%s4.2 #&#
\subsection{Iterative screening procedure}\label{sec4.2}
As we can see from the main results, the proposed marginal empirical
likelihood screening procedure works ideally for the case with
explanatory variables that are independent of each other. To deal
with challenging situations with correlated explanatory variables,
we propose to use the following iterative sure independence
screening procedure.

Step 1: %we use Equations (\ref{eqel}) and (\ref{eqelr}) to calculate
%marginal profile empirical likelihood $\mbox{EL}_j(\mu)$ and empirical
%likelihood ratio $\ell_j(\mu)$ for each predictor $j=1, 2, \ldots, p$.
Rank explanatory variables according to $\ell_j(0)$ by (\ref{eqellj})
%%(\textcolor{red}{giving same ranking order?})
and select top ranked explanatory variables with largest values of
$\ell_j(0)$'s until some desirable number of features are included.
Denote the set of
select explanatory variables by~$\widehat{\cal M}_1$.

Step 1$'$:\vspace*{1pt} Apply penalized empirical likelihood
[\citet{TangLeng2010Bioka,LengTang2010}] to explanatory variables in
$\widehat{\cal
M}_1$ and denote the final model by $\widehat{\cal A}_1$.

Step 2: Let $\widehat{\cal A}_k\subset\{1,\ldots,p\}$ be the
selected model at the $k$th step. At the $k$th iteration, for each
$j\notin\widehat{\cal A}_k $, denote by
\[
\mathrm{EL}_{\{j\}\cup
\widehat{\cal A}_k }(\bmu) =\sup\Biggl\{ \prod
_{i=1}^n w_i\dvtx w_i\ge0,
\sum_{i=1}^n w_i=1, \sum
_{i=1}^n w_i
\bX_{i,\{j\}\cup\widehat
{\cal A}_k }Y_i=\bmu\Biggr\}
\]
the empirical likelihood for the combined covariates, and denote by
\[
\widetilde{\mathrm{EL}}_{j}(\mu)=\sup_{\mu_j=\mu} \bigl\{
\mathrm{EL}_{\{
j\}\cup\widehat{\cal A}_k }(\bmu)\bigr\}
\]
the profile empirical likelihood evaluated at $\mu$.
%of predictor $j$ given predictors in $\widehat{\cal A}_1$ for any
%predictor $j\notin\widehat{\cal A}_k$.
Rank explanatory variable $j$ in $\widehat{\cal A}_k^c$ according to
$\widetilde{\mathrm{EL}}_{j}(0)$ and select the top ranked until some
desirable number of features are included. Denote the
selected set by $\widehat{\cal M}_{k+1}$.

Step 2$'$: Apply penalized empirical likelihood to explanatory
variables in $\widehat{\cal A}_k\cup\widehat{\cal M}_{k+1}$ and
denote the final model by $\widehat{\cal A}_{k+1}$.

Step 3: Repeat steps 2 and 2$'$ when either $\widehat{\cal
A}_{k+1}=\widehat{\cal A}_{k}$ or the size of $\widehat{\cal
A}_{k+1}$ reaches a pre-specified number.

The above iterative screening procedure incorporates the profile
empirical likelihood. The rationale behind it is to capture the
joint impact that may be invisible using the marginal screening
procedure if correlations exist among those covariates. Our
iterative screening procedure shares some similar features of the
analogous ones in \citet{FanLv2008JRSSB} and \citet{FanSong2010AOS}.
However, on the other hand, the iterative procedure using the
profile empirical likelihood ratio shares the feature of the
marginal empirical likelihood approach by incorporating the level of
uncertainties. In addition, we note that the above iterative
procedure is generally applicable in a broad framework. %to both linear
%models
%and the generalized linear model.

%s5 #&#
\section{Numerical examples}\label{s5}

In this section, we use five simulation examples and a real data
example to demonstrate the performance of the proposed empirical
likelihood-based screening procedure (denoted by EL-SIS) and
corresponding iterative procedure (denoted by EL-ISIS). Depending on
the example setting, we compare it with the screening methods
proposed in \citet{FanLv2008JRSSB} (denoted by LS-SIS and LS-ISIS) and
\citet{FanSong2010AOS} (denoted by GLM-SIS and GLM-ISIS) for linear
regression models and generalized linear models, respectively. Whenever
appropriate, we compare to the robust rank correlation based screening
(RRC-SIS and RRC-ISIS) studied by \citet{Lietal2012AOS}.
For
all simulation examples, we begin with $p=1000$ explanatory
variables and screen to a much smaller number $d$ of explanatory
variables. The respective SCAD penalized variable selection is
further applied to these selected explanatory variables to get the
corresponding final model. Results over 200 repetitions are
reported. For each case, we report the number of repetitions that
each important explanatory variable is selected in the final model
and also the average number of unimportant explanatory variables
being selected.

%ex1 #&#
%
\begin{example}\label{exampl1}
This example has a very standard setting with
three important explanatory variables and is taken from \citet
{FanLv2008JRSSB}. Covariates are generated as $X_j\sim N(0,1)$ and
$\operatorname{cov}(X_j, X_{j'})=1$ if $j=j'$ and $0.3$ otherwise. The
response is generated as $Y=5X_1+5X_2+5X_3+\varepsilon$ with error
being independent of the explanatory variables. We consider three
different error distribution $N(0,1)$, $N(0,2^2)$, and $t_4$ for
$\varepsilon$. Random samples of size $n=100$ are used and we set
$d=\lfloor n/(2\log n)\rfloor=10$, where $\lfloor a \rfloor$ denotes
the largest integer that is less than or equal to $a$. Results over
200 repetitions are reported in Table~\ref{tbex1}, where we report
%
%d=floor(n/log(n)/2)
%
%t1 #&#
%
\begin{table}
\tablewidth=280pt
%/home/statlab/yichaowu/TangCL/rerunEx2center/summary.R
\caption{Simulation result for Example \protect\ref{exampl1}}\label{tbex1}
\begin{tabular*}{\tablewidth}{@{\extracolsep{\fill}}lccccc@{}}
\hline
&&&&& \multicolumn{1}{c@{}}{\textbf{Unimportant explanatory}}\\
$\bolds{\varepsilon}$ & \multicolumn{1}{c}{\textbf{Method}}
& \multicolumn{1}{c}{$\bolds{X_1}$} & \multicolumn{1}{c}{$\bolds{X_2}$}
& \multicolumn{1}{c}{$\bolds{X_3}$} & \multicolumn{1}{c@{}}{\textbf{variables}} \\
\hline
$N(0,1)$ & LS-SIS & 199 & 199 & 200 &1.406219 \\
& RRC-SIS & 199 & 199 & 199 &1.407222\\
& EL-SIS & 194 & 183 & 185 &1.442327 \\
[3pt]
& LS-ISIS & 200 & 200 & 200 &0.965898 \\
& RRC-ISIS & 200 & 200 & 200 &0.800401\\
& EL-ISIS & 200 & 200 & 200 &0.659980 \\
[6pt]
$N(0,2^2)$ & LS-SIS & 199 & 199 & 200& 1.406219\\
& RRC-SIS & 198 & 198 & 199 & 1.409228\\
& EL-SIS & 192 & 182 & 183& 1.447342 \\
[3pt]
& LS-ISIS & 200 & 200 & 200& 1.404213 \\
& RRC-ISIS & 200 & 200 & 200 & 1.403210\\
& EL-ISIS & 200 & 200 & 200& 0.980943\\
[6pt]
$t_4$ & LS-SIS & 199 & 199 & 200 &1.406219\\
& RRC-SIS & 198 & 199 & 199 &1.408225\\
& EL-SIS & 193 & 186 & 187 &1.438315 \\
[3pt]
& LS-ISIS & 200 & 200 & 200 &1.383149\\
& RRC-ISIS & 200 & 200 & 200 &1.362086\\
& EL-ISIS & 200 & 199 & 200 &0.635908 \\
\hline
\end{tabular*}
\end{table}
the number of repetitions that each of the important explanatory
variables $X_1$, $X_2$ and $X_3$ is selected. For unimportant
explanatory variables, Table~\ref{tbex1} reports their average number
of repetitions for each being selected.
% and also the average
%number of repetitions other unimportant explanatory variables being
%selected.
It
shows that the proposed empirical likelihood-based screening methods
perform very competitively when compared to the least squares-based
screening or the robust rank correlation-based screening.
\end{example}

%ex2 #&#
%
\begin{example}\label{exampl2}
The second example is also from \citet{FanLv2008JRSSB}
and has a hidden important explanatory variable, which is important
but marginally uncorrelated with the response. This example is to
illustrate that the proposed iterative empirical likelihood-based
screening works effectively in such challenging cases. Covariates
are generated as $X_j\sim N(0,1)$ and $\operatorname{cov}(X_j,
X_{j'})=1$ if
$j=j'$ and $0.3$ otherwise except $\operatorname{cov}(X_4, X_j)=\sqrt{0.3}$
for $j\ne4$. The response is generated as
$Y=5X_1+5X_2+5X_3-15\sqrt{0.3}X_4+\varepsilon$ with $\varepsilon$
being independent of explanatory variables. We consider three
different error distribution $N(0,1)$, $N(0,2^2)$, and $t_4$.
Results over 200 repetitions with $n=100$ and $d=\lfloor n/(2\log
n)\rfloor=10$ are reported in Table~\ref{tbex2}. It shows that the
%
%d=floor(n/log(n)/2)
%
%t2 #&#
%
\begin{table}
%/home/statlab/yichaowu/TangCL/rerunEx2center/summary.R
\caption{Simulation result for Example \protect\ref{exampl2} with a
hidden important
explanatory variable $X_4$}\label{tbex2}
\begin{tabular*}{\tablewidth}{@{\extracolsep{\fill}}lccccd{3.0}c@{}}
\hline
&&&&&& \multicolumn{1}{c@{}}{\textbf{Unimportant explanatory}}\\
$\bolds{\varepsilon}$ & \textbf{Method} & $\bolds{X_1}$ & $\bolds{X_2}$
& $\bolds{X_3}$ & \multicolumn{1}{c}{$\bolds{X_4}$ \textbf{(hidden)}} &
\multicolumn{1}{c@{}}{\textbf{variables}} \\
\hline
$N(0,1)$ & LS-SIS & 198 & 197& 195& 0 &1.415663 \\
& RRC-SIS & 196 & 197& 194& 0 &1.418675\\
& EL-SIS & 198 & 198& 197& 0 &1.412651 \\ [3pt]
& LS-ISIS & 200 & 200& 199& 196 &1.125502 \\
& RRC-ISIS & 200 & 199& 200& 111 &1.157631\\
& EL-ISIS & 199 & 199& 200& 193 &0.853414 \\
[6pt]
$N(0,2^2)$ & LS-SIS & 198 & 197& 194& 0 &1.416667 \\
& RRC-SIS & 196 & 196& 194& 0 &1.419679\\
& EL-SIS & 198 & 196& 194& 0 &1.417671 \\ [3pt]
& LS-ISIS & 199 & 200& 199& 196 &1.210843 \\
& RRC-ISIS & 199 & 199 & 200 & 96 &1.311245\\
& EL-ISIS & 200 & 200& 198& 188 &0.912651 \\
[6pt]
$t_4$ & LS-SIS & 197 & 197& 197 & 0 &1.414659 \\
& RRC-SIS & 195 & 198 & 196 & 0 &1.416667\\
& EL-SIS & 197 & 198& 196 & 0 &1.414659 \\ [3pt]
& LS-ISIS & 199 & 200& 200 & 196 &1.209839 \\
& RRC-ISIS & 200& 200& 199& 100 &1.305221\\
& EL-ISIS & 200 & 198& 200 & 185 &0.824297 \\
\hline
\end{tabular*}
\end{table}
empirical likelihood-based screening is challenged by the hidden
important explanatory variable $X_4$ but the corresponding
iterative screening can easily pick it up. Overall the performance
of the empirical likelihood-based screening methods is very similar
to that of the least squares-based screening methods and is better than
the robust rank correlation-based screening. Note that iterative
version of the robust rank correlation-based screening is
residual-based. This explains the improvement of the robust rank
correlation-based screening.
\end{example}

%ex3 #&#
%
\begin{example}\label{exampl3}
The performances of the empirical
likelihood-based screening and the least squares-based screening
methods are very similar in the previous two examples. It is known
that the empirical likelihood approach requires a less restrictive
distributional assumption. We next use a heteroscedastic example to
show the advantage of the empirical likelihood-based screening.
Explanatory variables are generated as $X_j\sim N(0,1)$ with
$\operatorname{cov}(X_j, X_{j'})=0$ for $j\ne j'$. The response is generated
as $Y=c(X_1-X_2+X_3)+\varepsilon/(X_1^2+X_2^2+X_3^2)$ with
independent $\varepsilon\sim N(0,1)$ and $c>0$ controls the signal
level. Results over 200 repetitions with $n=70$ and $d=\lfloor
n/(2\log n)\rfloor=8$ are reported in Table~\ref{tbex3} for three
%
%d=floor(n/log(n)/2)
%/home/statlab/yichaowu/TangCL/twogoodex
%
%t3 #&#
%
\begin{table}
\tablewidth=280pt
\caption{Simulation result for Example \protect\ref{exampl3}}\label{tbex3}
\begin{tabular*}{\tablewidth}{@{\extracolsep{\fill}}lccccc@{}}
\hline
&&&&& \multicolumn{1}{c@{}}{\textbf{Unimportant explanatory}}\\
$\bolds{c}$ & \textbf{Method} & $\bolds{X_1}$ & $\bolds{X_2}$ &
$\bolds{X_3}$ & \textbf{variables} \\
\hline
1 & LS-SIS &149 & 147 & 156 &1.151454 \\
& RRC-SIS & 191 & 185& 190& 1.037111\\
& EL-SIS &190 & 184 &191 &1.038114 \\
% &aEL-SIS &191& 194& 195& 2.18656\\
[4pt]
1.5 & LS-SIS & 173 & 171& 174 &1.085256\\
& RRC-SIS & 194 & 191& 193& 1.025075\\
& EL-SIS & 196 & 192 & 194& 1.021063 \\
% & aEL-SIS &194& 198& 198& 2.15346\\
[4pt]
2 & LS-SIS & 182 & 182 & 180& 1.059178 \\
& RRC-SIS & 194 & 194 & 195& 1.020060\\
& EL-SIS & 199 & 195 & 194 &1.015045 \\
% &aEL-SIS &198& 200& 199& 2.171515\\
\hline
\end{tabular*}
\end{table}
different values of $c$. It shows that the performance of the least
squares-based screening is severely affected by the
heteroscedasticity especially when the signal level is low. On the
other hand, the proposed empirical likelihood-based screening works
much better and similarly as the robust rank correlation-based screening.
\end{example}

%ex4 #&#
%
\begin{example}\label{exampl4}
%Can you write down the details how the data are generated? For your
%reference, here I copy and paste the code for data generate
%(commented). Results over 200 repetitions with $n=60$ and $d=15$ are
%reported in Table \ref{tbex4}.}
We now consider an example with the extended scope. In this example, we
generate data from the longitudinal data example as in Section~\ref
{sec4.2} with $m=4$ means 4 repeated measurements generated. In
particular, the following model is generated:
\[
Y_{il}=\bX_{il}^{\mathrm{T}}\bbeta+\varepsilon_{il}
\qquad(i=1,\ldots,n; l=1,\ldots,m).
\]
Here $\bX_{il}$ is generated from multivariate normal
$N(\bzero,{\bSigma})$ with $\bSigma=(\sigma_{jk})_{j,k=1,\ldots,p}$
and $\sigma_{jk}=0.5^{|j-k|}$. The error vector
$\bvare_i=(\varepsilon_{i1},\ldots,\varepsilon_{im})^{\mathrm{T}}$
is generated
from multivariate normal distribution with unit variance.
The correlation structure of $\bvare$ is specified as AR(1) with
parameter 0.8; see \citet{Diggleetal2002}
for reference for the correlation structure.
The first five components of the true $\bbeta$ is set to be $c\cdot
(2.0,-2.0,0,0,2.0)^{\mathrm{T}}$
where $c$ is used to control the signal strength, and all other
components of $\bbeta$ are zero.
We use two sets of basis matrices in (\ref{eqgee2}). We
take $\bM_1=\bI$ as the identity matrix. The second basis matrix
$\bM_2$ is a matrix with two main off-diagonals being 1 and 0
elsewhere corresponding to the AR(1) working correlation
[\citet{Quetal2000Bioka}]. We then apply the marginal empirical
likelihood procedure as
in Section~\ref{s4} using the marginal estimating function of
(\ref{eqgee2}). Here we note that the marginal estimating function
is 4-dimensional. By ignoring the correlation structure of the
longitudinal data, the least squares-based screening and robust rank
correlation-based screening procedures can
be applied. Results over 200 repetitions with $n=60$ and $d=15$ are
reported in Table~\ref{tbex4}. From Table~\ref{tbex4}, we clearly
see that the marginal empirical likelihood approach works much
better than the alternative ones, especially when signal is
relatively weak. The improvement can be seen as the results of
incorporating additional data structural information. Hence, we
demonstrate an advantage of the marginal empirical likelihood
approach of being adaptive and flexible.

%
%t4 #&#
%
\begin{table}
\tablewidth=280pt
%/home/statlab/yichaowu/TangCL/twogoodex
\caption{Simulation result for the longitudinal data estimation
function example with $c$ controlling the signal
strength}\label{tbex4}
\begin{tabular*}{\tablewidth}{@{\extracolsep{\fill}}lccccc@{}}
\hline
&&&&& \multicolumn{1}{c@{}}{\textbf{Unimportant explanatory}}\\
$\bolds{c}$ & \textbf{Method} & $\bolds{X_1}$ & $\bolds{X_2}$ &
$\bolds{X_5}$ & \textbf{variables} \\
\hline
1 & LS-SIS & \hphantom{0}90 & \hphantom{0}73 & 153 & 2.692076\\
& GEE-SIS & 111 & 111 & 168 &2.617854\\
& RRC-SIS & \hphantom{0}84 & \hphantom{0}66 & 136& 2.722166\\
& EL-SIS & 135 & 128& 191& 2.553661\\
% & aEL-SIS & 136& 128& 194& 2.451354\\
[4pt]
1.5 & LS-SIS & 153 & 153& 195& 2.506520\\
&GEE-SIS &165& 160 & 196 & 2.486459\\
& RRC-SIS & 142 & 136 & 193& 2.536610\\
& EL-SIS & 176 & 187& 199& 2.445336\\
% &aEL-SIS &176& 177& 200& 2.381143\\
[4pt]
2 & LS-SIS & 183 & 183& 200& 2.441324\\
&GEE-SIS & 183& 184& 200& 2.440321\\
& RRC-SIS & 179 & 176 & 200 &2.452357\\
& EL-SIS & 192 & 196& 200& 2.419258\\
% &aEL-SIS &194& 192& 200& 2.382146\\
[4pt]
2.5 & LS-SIS & 195 & 195& 200& 2.417252\\
&GEE-SIS & 196& 195& 200& 2.416249\\
& RRC-SIS & 192 & 190 & 200 &2.425276\\
& EL-SIS & 198 & 197& 200& 2.412237\\
% &aEL-SIS &199& 196& 200& 2.363089\\
[4pt]
3 & LS-SIS & 199 & 198& 200& 2.410231\\
&GEE-SIS &198& 197& 200& 2.412237\\
& RRC-SIS & 199 & 198 & 200& 2.410231\\
& EL-SIS & 200 & 198& 200& 2.409228\\
% &aEL-SIS & 199& 200& 200& 2.399198\\
\hline
\end{tabular*}
\end{table}

In the review process, one referee pointed out that our comparison to
the LS-SIS is not fair as it is based on the ordinary least squares. It
is more reasonable to compare to a weighted least squares-based
screening by adjusting to correlation among longitudinal observations.
To address this issue, we implement this weighted least squares-based
screening by using the R package ``geepack,'' which can estimate both
the correlation structure and regression parameter once a parametric
form of the correlation structure is specified. Table~\ref{tbex4} is
updated accordingly with GEE-SIS denoting this weighted least
squares-based screening method. It shows that our newly proposed EL-SIS
still does better than the GEE-SIS even though a correct parametric
correlation structure, AR(1), is specified.
\end{example}

%d=floor(sqrt(n*ni))

%ni=4 # the number of observations for each subject
%n=60 # sample size
%p=1000 # dimensionality
%r=2*p # number of estimating equations
%pbeta=2*c(1.0,-1.0,0,0,1.0,rep(0,p-5)) # true beta
%
%nrep=200
%
%d=floor(sqrt(n*ni))
%
%rho=0.3 # the correlation in generating X
%a=0.8 # the correlation in generating e
%
%## Matrices in Qu et al. (2000)
%cor.type='ar' ## exchangeable
%MM=mtrx(ni,type=cor.type)
%M1=MM$M1
%M2=MM$M2
%
%### The 2p EEs
%m=function(theta,xy) {
% p=ncol(xy)-1
% x=xy[,2:(p+1)]
% y=xy[,1]
% r=2*p
% z=matrix(0,n,r)
% idx=1:ni
% for(i in 1:n) {
% if(p>1){
% xi=x[idx,]
% res=y[idx]-xi%*%theta
% } else
% {xi=x[idx]
% res=y[idx]-xi*theta
% }
% z[i,1:p]= drop(t(xi)%*%M1%*%res)
% z[i,(p+1):(2*p)]=drop(t(xi)%*%M2%*%res)
% idx=idx+ni
% }
% z
%}
%
%lsSIS=matrix(0,p,nrep)
%elSIS=matrix(0,p,nrep)
%
%vm=diag(rep(4.5,ni))
%
%## DGP longitudinal reponses
%
%for (irep in 1:nrep) {
%cat(irep,'\n')
%
%x=mvrnorm(n*ni,rep(0,p),ar1(rho,p))
%e=rep(0,n*ni)
%idx=1:ni
%for(j in 1:n) {
%e[idx]=mvrnorm(1,rep(0,ni),vm%*%ar1(a,ni)%*%vm)
%idx=idx+ni
%}
%y=x%*%pbeta+e

%ex5 #&#
%
\begin{example}\label{exampl5}
This is an extension of Example \ref{exampl2} to the case
with a binary response using logistic regression. Covariates are
generated as $X_j\sim N(0,1)$ and $\operatorname{cov}(X_j, X_{j'})=1$ if
$j=j'$ and $0.3$ otherwise except $\operatorname{cov}(X_4, X_j)=\sqrt{0.3}$
for $j\ne4$. The binary response is generated from Bernoulli
distribution with success probability given by
$\{1+\exp(-4X_1-4X_2-4X_3+12\sqrt{0.3}X_4)\}^{-1}$. Results over 200
repetitions with $n=400$ and $d=10$ are reported in Table~\ref{tbex5}.
%
%t5 #&#
%
\begin{table}
%%%%%%/home/statlab/yichaowu/TangCL/logisticCOR3centerRangeNewRange/trynewrange
\caption{Simulation result for Example \protect\ref{exampl5}}\label{tbex5}
\begin{tabular*}{\tablewidth}{@{\extracolsep{\fill}}lcccd{3.0}c@{}}
\hline
&&&&& \multicolumn{1}{c@{}}{\textbf{Unimportant explanatory}}\\
\textbf{Method} & $\bolds{X_1}$ & $\bolds{X_2}$
& $\bolds{X_3}$ & \multicolumn{1}{c}{$\bolds{X_4}$ \textbf{(hidden)}} &
\multicolumn{1}{c@{}}{\textbf{variables}} \\
\hline
GLM-SIS & 200 & 200 &200 & 0& 1.405623 \\
RRC-SIS & 200& 200& 200& 0& 1.405623\\
EL-SIS & 200 & 200 &200 & 0& 1.405623 \\[4pt]
GLM-ISIS & 200 & 200 &200 & 200& 0.324297 \\
EL-ISIS & 199 & 200 &200 & 199& 0.764056 \\
\hline
\end{tabular*}
\end{table}
A similar performance pattern is observed. For this
example, the result for the iterative version of the robust rank
correlation-based screening is not presented since it is not clear how
to define a residual-based iterative procedure.
\end{example}

\textit{A real data example}.
%/home/statlab/yichaowu/TangCL/RealGSE5680/realdata2outcenter.Rdata
Glioblastoma is the most common primary malignant brain tumor of
adults and one of the most lethal of all cancers
[\citet{Horvathetal2006PNAS}]. The median survival of glioblastoma
patients is 15 months
from the time of diagnosis. We next apply our proposed methods to a
microarray gene expression dataset of glioblastoma patients reported
in \citet{Horvathetal2006PNAS}. The dataset has been analyzed by
\citet
{Pan2010Biometrics} and \citet{LiLi2008Biominformatics} among many
others. Drawn from two
different studies, the data consist of two independent sets. We use
the set with 50 samples. We use the log survival time, measured in
years, as the response. The second sample with a outlier response
is excluded and the other 49 samples are used in our analysis.
Explanatory variables are gene expression profiles of 1523 genes
measured on Affymetrix HG-U133A arrays.

We apply the least squares-based and empirical likelihood-based
screening methods with $d=6$. LS-SIS selects ``GSN'', ``FOS'',
``COL11A1'', ``AVPR1A'', ``SELE'', and ``TBL1X'' as important gene
explanatory variables while EL-SIS selects ``GSN'', ``JAK2'',
``COL11A1'', ``CDK6'', ``ADCYAP1R1'', and ``TBL1X''. Note that
they select some common genes (``GSN'' and\break ``COL11A1'') and some
different genes. LS-ISIS selects ``GSN'', ``COL11A1'', ``THBS1'',
``SELE'', ``TBL1X'', and ``GCGR''. EL-ISIS selects ``DUSP7'',
``COL11A1'', ``BST1'', ``ADCYAP1R1'', ``TBL1X'', and ``GCGR''.
Similarly two genes (``TBL1X'' and ``GCGR'') are recruited by the
iterative screening methods based on both the least squares and
empirical likelihood. The robust rank correlation-based screening
performs similarly with 2--3 overlapping genes.

%LS-SIS selects 118 557 844 956 1134 1382
%LS-ISIS selects 118 844 1040 1134 1382 1489
%EL-SIS selects 118 326 844 845 1352 1382
%EL-ISIS selects 120 844 1329 1352 1382 1489
%p=1523 n=50-1

%> load("realdata2outcenterRank.Rdata")
%> which(lsSIS[,1]>0.5)
%[1] 42 118 550 1192 1352 1382
%> which(lsISIS[,1]>0.5)
%[1] 42 120 232 476 1352 1382

%s6 #&#
\section{Discussion}\label{s6}

Screening based on marginal model fitting has enjoyed great
popularity in the recent literature. However, most, if not all, of
the marginal screening methods studied thus far are based on some
restrictive distributional assumptions. Yet these assumptions may
not be realistic in applications. Thus motivated we propose a new
screening method based on marginal empirical likelihood, which is
known to be less restrictive. It has been demonstrated to be
effective using both theoretical sure screening property and
numerical evidences. Further extensions using empirical likelihood
are being investigated.

\section*{Acknowledgments}

We thank the Editor, the Associate Editor and two referees for very
constructive comments and suggestions which have improved the
presentation of the paper. We are very grateful to Drs. Gaorong Li, Yi
Li, Heng Peng and Sihai Dave Zhao for sharing with us programs for
implementing their methods and Drs. Hongzhe Li and Wei Pan for sharing
the real data.

\begin{supplement}%[id=suppA]
\stitle{Supplement to ``Marginal empirical
likelihood and sure independence feature screening''}
\slink[doi]{10.1214/13-AOS1139SUPP} %[doi,text={...}] - jei reikia
%suskaldyti doi
\sdatatype{.pdf}
\sfilename{aos1139\_supp.pdf}
\sdescription{This supplement contains all technical proofs.}
\end{supplement}

% imsref loaded by lrinkeviciute, 2013-08-20 08:36:52
%
% imsref loaded by lrinkeviciute, 2013-08-20 10:05:10
%

\printaddresses


\begin{thebibliography}{38}
% BibTex style file: ims.bst, 2013-01-28
% Default style options (sort=0,type=number).
% Used options (sort=1,type=nameyear).

%b1 #&#
\bibitem[\protect\citeauthoryear{B{\"u}hlmann and van~de
Geer}{2011}]{Buhlmannvan2011}
%
\begin{bbook}[mr]
\bauthor{\bsnm{B{\"u}hlmann},~\bfnm{Peter}\binits{P.}} \AND
\bauthor{\bparticle{van~de} \bsnm{Geer},~\bfnm{Sara}\binits{S.}}
(\byear{2011}).
\btitle{Statistics for High-dimensional Data:
Methods, Theory and Applications}.
\bpublisher{Springer}, \blocation{Heidelberg}.
\bid{doi={10.1007/978-3-642-20192-9}, mr={2807761}}
\bptok{imsref}%
\end{bbook}
%
\endbibitem

%b2 #&#
\bibitem[\protect\citeauthoryear{Chang, Chen and Chen}{2013}]{Changetal2012}
%
\begin{bmisc}[auto:STB|2013/06/05|13:45:01]
\bauthor{\bsnm{Chang},~\bfnm{J.}\binits{J.}},
\bauthor{\bsnm{Chen},~\bfnm{S.~X.}\binits{S.~X.}} \AND
\bauthor{\bsnm{Chen},~\bfnm{X.}\binits{X.}}
(\byear{2013}).
\bhowpublished{High dimensional generalized empirical likelihood for moment
restrictions with dependent data.
Available at \arxivurl{arXiv:1308.5732}}.
\bptok{imsref}%
\end{bmisc}
%
\endbibitem

%b3 #&#
\bibitem[\protect\citeauthoryear{Chang, Tang and Wu}{2013}]{Changetal2013}
%
\begin{bmisc}[auto:STB|2013/06/05|13:45:01]
\bauthor{\bsnm{Chang},~\bfnm{J.}\binits{J.}},
\bauthor{\bsnm{Tang},~\bfnm{C.~Y.}\binits{C.~Y.}} \AND
\bauthor{\bsnm{Wu},~\bfnm{Y.}\binits{Y.}}
(\byear{2013}).
\bhowpublished{Supplement to ``Marginal empirical likelihood and sure
independence feature screening.'' DOI:\doiurl{10.1214/13-AOS1139SUPP}}.
\bptok{imsref}%
\end{bmisc}
%
\endbibitem

%b4 #&#
\bibitem[\protect\citeauthoryear{Chen and Cui}{2003}]{ChenCui2003Sinica}
%
\begin{barticle}[mr]
\bauthor{\bsnm{Chen},~\bfnm{Song~Xi}\binits{S.~X.}} \AND
\bauthor{\bsnm{Cui},~\bfnm{Hengjian}\binits{H.}}
(\byear{2003}).
\btitle{An extended empirical likelihood for generalized linear models}.
\bjournal{Statist. Sinica}
\bvolume{13}
\bpages{69--81}.
\bid{issn={1017-0405}, mr={1963920}}
\bptok{imsref}%
\end{barticle}
%
\endbibitem

%b5 #&#
\bibitem[\protect\citeauthoryear{Chen, Gao and Tang}{2008}]{ChenGaoTang2008}
%
\begin{barticle}[mr]
\bauthor{\bsnm{Chen},~\bfnm{Song~Xi}\binits{S.~X.}},
\bauthor{\bsnm{Gao},~\bfnm{Jiti}\binits{J.}} \AND
\bauthor{\bsnm{Tang},~\bfnm{Cheng~Yong}\binits{C.~Y.}}
(\byear{2008}).
\btitle{A test for model specification of diffusion processes}.
\bjournal{Ann. Statist.}
\bvolume{36}
\bpages{167--198}.
\bid{doi={10.1214/009053607000000659}, issn={0090-5364}, mr={2387968}}
\bptok{imsref}%
\end{barticle}
%
\endbibitem

%b6 #&#
\bibitem[\protect\citeauthoryear{Chen, Peng and Qin}{2009}]{ChenPengQin2008}
%
\begin{barticle}[mr]
\bauthor{\bsnm{Chen},~\bfnm{Song~Xi}\binits{S.~X.}},
\bauthor{\bsnm{Peng},~\bfnm{Liang}\binits{L.}} \AND
\bauthor{\bsnm{Qin},~\bfnm{Ying-Li}\binits{Y.-L.}}
(\byear{2009}).
\btitle{Effects of data dimension on empirical likelihood}.
\bjournal{Biometrika}
\bvolume{96}
\bpages{711--722}.
\bid{doi={10.1093/biomet/asp037}, issn={0006-3444}, mr={2538767}}
\bptok{imsref}%
\end{barticle}
%
\endbibitem

%b7 #&#
\bibitem[\protect\citeauthoryear{Chen and Van~Keilegom}{2009}]{ChenVan2009test}
%
\begin{barticle}[mr]
\bauthor{\bsnm{Chen},~\bfnm{Song~Xi}\binits{S.~X.}} \AND
\bauthor{\bsnm{Van~Keilegom},~\bfnm{Ingrid}\binits{I.}}
(\byear{2009}).
\btitle{A review on empirical likelihood methods for regression}.
\bjournal{TEST}
\bvolume{18}
\bpages{415--447}.
\bid{doi={10.1007/s11749-009-0159-5}, issn={1133-0686}, mr={2566404}}
\bptok{imsref}%
\end{barticle}
%
\endbibitem

%b8 #&#
\bibitem[\protect\citeauthoryear{Diggle et~al.}{2002}]{Diggleetal2002}
%
\begin{bbook}[mr]
\bauthor{\bsnm{Diggle},~\bfnm{Peter~J.}\binits{P.~J.}},
\bauthor{\bsnm{Heagerty},~\bfnm{Patrick~J.}\binits{P.~J.}},
\bauthor{\bsnm{Liang},~\bfnm{Kung-Yee}\binits{K.-Y.}} \AND
\bauthor{\bsnm{Zeger},~\bfnm{Scott~L.}\binits{S.~L.}}
(\byear{2002}).
\btitle{Analysis of Longitudinal Data},
\bedition{2nd} ed.
\bseries{Oxford Statistical Science Series}
\bvolume{25}.
\bpublisher{Oxford Univ. Press}, \blocation{Oxford}.
\bid{mr={2049007}}
\bptok{imsref}%
\end{bbook}
%
\endbibitem

%b9 #&#
\bibitem[\protect\citeauthoryear{Fan, Feng and
Song}{2011}]{FangFengSongJASA2011}
%
\begin{barticle}[mr]
\bauthor{\bsnm{Fan},~\bfnm{Jianqing}\binits{J.}},
\bauthor{\bsnm{Feng},~\bfnm{Yang}\binits{Y.}} \AND
\bauthor{\bsnm{Song},~\bfnm{Rui}\binits{R.}}
(\byear{2011}).
\btitle{Nonparametric independence screening in sparse ultra-high-dimensional
additive models}.
\bjournal{J. Amer. Statist. Assoc.}
\bvolume{106}
\bpages{544--557}.
\bid{doi={10.1198/jasa.2011.tm09779}, issn={0162-1459}, mr={2847969}}
\bptok{imsref}%
\end{barticle}
%
\endbibitem

%b10 #&#
\bibitem[\protect\citeauthoryear{Fan and Lv}{2008}]{FanLv2008JRSSB}
%
\begin{barticle}[mr]
\bauthor{\bsnm{Fan},~\bfnm{Jianqing}\binits{J.}} \AND
\bauthor{\bsnm{Lv},~\bfnm{Jinchi}\binits{J.}}
(\byear{2008}).
\btitle{Sure independence screening for ultrahigh dimensional feature space}.
\bjournal{J.~R. Stat. Soc. Ser. B Stat. Methodol.}
\bvolume{70}
\bpages{849--911}.
\bid{doi={10.1111/j.1467-9868.2008.00674.x}, issn={1369-7412}, mr={2530322}}
\bptok{imsref}%
\end{barticle}
%
\endbibitem

%b11 #&#
\bibitem[\protect\citeauthoryear{Fan and Lv}{2010}]{FanLv2009Sinica}
%
\begin{barticle}[mr]
\bauthor{\bsnm{Fan},~\bfnm{Jianqing}\binits{J.}} \AND
\bauthor{\bsnm{Lv},~\bfnm{Jinchi}\binits{J.}}
(\byear{2010}).
\btitle{A selective overview of variable selection in high dimensional feature
space}.
\bjournal{Statist. Sinica}
\bvolume{20}
\bpages{101--148}.
\bid{issn={1017-0405}, mr={2640659}}
\bptok{imsref}%
\end{barticle}
%
\endbibitem

%b12 #&#
\bibitem[\protect\citeauthoryear{Fan and Song}{2010}]{FanSong2010AOS}
%
\begin{barticle}[mr]
\bauthor{\bsnm{Fan},~\bfnm{Jianqing}\binits{J.}} \AND
\bauthor{\bsnm{Song},~\bfnm{Rui}\binits{R.}}
(\byear{2010}).
\btitle{Sure independence screening in generalized linear models with
{NP}-dimensionality}.
\bjournal{Ann. Statist.}
\bvolume{38}
\bpages{3567--3604}.
\bid{doi={10.1214/10-AOS798}, issn={0090-5364}, mr={2766861}}
\bptok{imsref}%
\end{barticle}
%
\endbibitem

%b13 #&#
\bibitem[\protect\citeauthoryear{Hansen}{1982}]{Hansen1982ecnoca}
%
\begin{barticle}[mr]
\bauthor{\bsnm{Hansen},~\bfnm{Lars~Peter}\binits{L.~P.}}
(\byear{1982}).
\btitle{Large sample properties of generalized method of moments estimators}.
\bjournal{Econometrica}
\bvolume{50}
\bpages{1029--1054}.
\bid{doi={10.2307/1912775}, issn={0012-9682}, mr={0666123}}
\bptok{imsref}%
\end{barticle}
%
\endbibitem

%b14 #&#
\bibitem[\protect\citeauthoryear{Hastie, Tibshirani and
Friedman}{2009}]{Hastieetal2009book}
%
\begin{bbook}[mr]
\bauthor{\bsnm{Hastie},~\bfnm{Trevor}\binits{T.}},
\bauthor{\bsnm{Tibshirani},~\bfnm{Robert}\binits{R.}} \AND
\bauthor{\bsnm{Friedman},~\bfnm{Jerome}\binits{J.}}
(\byear{2009}).
\btitle{The Elements of Statistical Learning:
Data Mining, Inference, and Prediction},
\bedition{2nd} ed.
\bpublisher{Springer}, \blocation{New York}.
\bid{doi={10.1007/978-0-387-84858-7}, mr={2722294}}
\bptok{imsref}%
\end{bbook}
%
\endbibitem

%b15 #&#
\bibitem[\protect\citeauthoryear{Hjort, McKeague and
Van~Keilegom}{2009}]{Hjortetal2008AS}
%
\begin{barticle}[mr]
\bauthor{\bsnm{Hjort},~\bfnm{Nils~Lid}\binits{N.~L.}},
\bauthor{\bsnm{McKeague},~\bfnm{Ian~W.}\binits{I.~W.}} \AND
\bauthor{\bsnm{Van~Keilegom},~\bfnm{Ingrid}\binits{I.}}
(\byear{2009}).
\btitle{Extending the scope of empirical likelihood}.
\bjournal{Ann. Statist.}
\bvolume{37}
\bpages{1079--1111}.
\bid{doi={10.1214/07-AOS555}, issn={0090-5364}, mr={2509068}}
\bptok{imsref}%
\end{barticle}
%
\endbibitem

%b16 #&#
\bibitem[\protect\citeauthoryear{Horvath et~al.}{2006}]{Horvathetal2006PNAS}
%
\begin{barticle}[auto:STB|2013/06/05|13:45:01]
\bauthor{\bsnm{Horvath},~\bfnm{S.}\binits{S.}},
\bauthor{\bsnm{Zhang},~\bfnm{B.}\binits{B.}},
\bauthor{\bsnm{Carlson},~\bfnm{M.}\binits{M.}},
\bauthor{\bsnm{Lu},~\bfnm{K.~V.}\binits{K.~V.}},
\bauthor{\bsnm{Zhu},~\bfnm{S.}\binits{S.}},
\bauthor{\bsnm{Felciano},~\bfnm{R.~M.}\binits{R.~M.}},
\bauthor{\bsnm{Laurance},~\bfnm{M.~F.}\binits{M.~F.}},
\bauthor{\bsnm{Zhao},~\bfnm{W.}\binits{W.}},
\bauthor{\bsnm{Shu},~\bfnm{Q.}\binits{Q.}},
\bauthor{\bsnm{Lee},~\bfnm{Y.}\binits{Y.}},
\bauthor{\bsnm{Scheck},~\bfnm{A.~C.}\binits{A.~C.}},
\bauthor{\bsnm{Liau},~\bfnm{L.~M.}\binits{L.~M.}},
\bauthor{\bsnm{Wu},~\bfnm{H.}\binits{H.}},
\bauthor{\bsnm{Geschwind},~\bfnm{D.~H.}\binits{D.~H.}},
\bauthor{\bsnm{Febbo},~\bfnm{P.~G.}\binits{P.~G.}},
\bauthor{\bsnm{Kornblum},~\bfnm{H.~I.}\binits{H.~I.}},
\bauthor{\bsnm{Cloughesy},~\bfnm{T.~F.}\binits{T.~F.}},
\bauthor{\bsnm{Nelson},~\bfnm{S.~F.}\binits{S.~F.}} \AND
\bauthor{\bsnm{Mischel},~\bfnm{P.~S.}\binits{P.~S.}}
(\byear{2006}).
\btitle{Analysis of oncogenic signaling networks in glioblastoma identifies
ASPM as a molecular target}.
\bjournal{Proc. Natl. Acad. Sci. USA}
\bvolume{103}
\bpages{17402--17407}.
\bptok{imsref}%
\end{barticle}
%
\endbibitem

%b17 #&#
\bibitem[\protect\citeauthoryear{Huang, Horowitz and
Ma}{2008}]{HuangHortowitzMa2008AOS}
%
\begin{barticle}[mr]
\bauthor{\bsnm{Huang},~\bfnm{Jian}\binits{J.}},
\bauthor{\bsnm{Horowitz},~\bfnm{Joel~L.}\binits{J.~L.}} \AND
\bauthor{\bsnm{Ma},~\bfnm{Shuangge}\binits{S.}}
(\byear{2008}).
\btitle{Asymptotic properties of bridge estimators in sparse high-dimensional
regression models}.
\bjournal{Ann. Statist.}
\bvolume{36}
\bpages{587--613}.
\bid{doi={10.1214/009053607000000875}, issn={0090-5364}, mr={2396808}}
\bptok{imsref}%
\end{barticle}
%
\endbibitem

%b18 #&#
\bibitem[\protect\citeauthoryear{Kolaczyk}{1994}]{Kolaczyk1994Sinica}
%
\begin{barticle}[mr]
\bauthor{\bsnm{Kolaczyk},~\bfnm{Eric~D.}\binits{E.~D.}}
(\byear{1994}).
\btitle{Empirical likelihood for generalized linear models}.
\bjournal{Statist. Sinica}
\bvolume{4}
\bpages{199--218}.
\bid{issn={1017-0405}, mr={1282871}}
\bptok{imsref}%
\end{barticle}
%
\endbibitem

%b19 #&#
\bibitem[\protect\citeauthoryear{Leng and Tang}{2012}]{LengTang2010}
%
\begin{barticle}[mr]
\bauthor{\bsnm{Leng},~\bfnm{Chenlei}\binits{C.}} \AND
\bauthor{\bsnm{Tang},~\bfnm{Cheng~Yong}\binits{C.~Y.}}
(\byear{2012}).
\btitle{Penalized empirical likelihood and growing dimensional general
estimating equations}.
\bjournal{Biometrika}
\bvolume{99}
\bpages{703--716}.
\bid{doi={10.1093/biomet/ass014}, issn={0006-3444}, mr={2966779}}
\bptok{imsref}%
\end{barticle}
%
\endbibitem

%b20 #&#
\bibitem[\protect\citeauthoryear{Li and Li}{2008}]{LiLi2008Biominformatics}
%
\begin{barticle}[auto:STB|2013/06/05|13:45:01]
\bauthor{\bsnm{Li},~\bfnm{C.}\binits{C.}} \AND
\bauthor{\bsnm{Li},~\bfnm{H.}\binits{H.}}
(\byear{2008}).
\btitle{Network-constrained regularization and variable selection for analysis
of genomic data}.
\bjournal{Bioinformatics}
\bvolume{24}
\bpages{1175--1182}.
\bptok{imsref}%
\end{barticle}
%
\endbibitem

%b21 #&#
\bibitem[\protect\citeauthoryear{Li, Zhong and Zhu}{2012}]{Lietal2012JASA}
%
\begin{barticle}[mr]
\bauthor{\bsnm{Li},~\bfnm{Runze}\binits{R.}},
\bauthor{\bsnm{Zhong},~\bfnm{Wei}\binits{W.}} \AND
\bauthor{\bsnm{Zhu},~\bfnm{Liping}\binits{L.}}
(\byear{2012}).
\btitle{Feature screening via distance correlation learning}.
\bjournal{J. Amer. Statist. Assoc.}
\bvolume{107}
\bpages{1129--1139}.
\bid{doi={10.1080/01621459.2012.695654}, issn={0162-1459}, mr={3010900}}
\bptok{imsref}%
\end{barticle}
%
\endbibitem

%b22 #&#
\bibitem[\protect\citeauthoryear{Li et~al.}{2012}]{Lietal2012AOS}
%
\begin{barticle}[mr]
\bauthor{\bsnm{Li},~\bfnm{Gaorong}\binits{G.}},
\bauthor{\bsnm{Peng},~\bfnm{Heng}\binits{H.}},
\bauthor{\bsnm{Zhang},~\bfnm{Jun}\binits{J.}} \AND
\bauthor{\bsnm{Zhu},~\bfnm{Lixing}\binits{L.}}
(\byear{2012}).
\btitle{Robust rank correlation based screening}.
\bjournal{Ann. Statist.}
\bvolume{40}
\bpages{1846--1877}.
\bid{doi={10.1214/12-AOS1024}, issn={0090-5364}, mr={3015046}}
\bptok{imsref}%
\end{barticle}
%
\endbibitem

%b23 #&#
\bibitem[\protect\citeauthoryear{Liang and Zeger}{1986}]{LiangZeger1986Bioka}
%
\begin{barticle}[mr]
\bauthor{\bsnm{Liang},~\bfnm{Kung~Yee}\binits{K.~Y.}} \AND
\bauthor{\bsnm{Zeger},~\bfnm{Scott~L.}\binits{S.~L.}}
(\byear{1986}).
\btitle{Longitudinal data analysis using generalized linear models}.
\bjournal{Biometrika}
\bvolume{73}
\bpages{13--22}.
\bid{doi={10.1093/biomet/73.1.13}, issn={0006-3444}, mr={0836430}}
\bptok{imsref}%
\end{barticle}
%
\endbibitem

%b24 #&#
\bibitem[\protect\citeauthoryear{McCullagh and Nelder}{1989}]{GLM1989MN}
%
\begin{bbook}[mr]
\bauthor{\bsnm{McCullagh},~\bfnm{P.}\binits{P.}} \AND
\bauthor{\bsnm{Nelder},~\bfnm{J.~A.}\binits{J.~A.}}
(\byear{1989}).
\btitle{Generalized Linear Models}.
\bpublisher{Chapman \& Hall/CRC}, \blocation{New York}.
\bptok{imsref}%
\end{bbook}
%
\endbibitem

%b25 #&#
\bibitem[\protect\citeauthoryear{Newey and Smith}{2004}]{NeweySmith2004Econca}
%
\begin{barticle}[mr]
\bauthor{\bsnm{Newey},~\bfnm{Whitney~K.}\binits{W.~K.}} \AND
\bauthor{\bsnm{Smith},~\bfnm{Richard~J.}\binits{R.~J.}}
(\byear{2004}).
\btitle{Higher order properties of {GMM} and generalized empirical likelihood
estimators}.
\bjournal{Econometrica}
\bvolume{72}
\bpages{219--255}.
\bid{doi={10.1111/j.1468-0262.2004.00482.x}, issn={0012-9682}, mr={2031017}}
\bptok{imsref}%
\end{barticle}
%
\endbibitem

%b26 #&#
\bibitem[\protect\citeauthoryear{Owen}{1988}]{Owen1988}
%
\begin{barticle}[mr]
\bauthor{\bsnm{Owen},~\bfnm{Art~B.}\binits{A.~B.}}
(\byear{1988}).
\btitle{Empirical likelihood ratio confidence intervals for a single
functional}.
\bjournal{Biometrika}
\bvolume{75}
\bpages{237--249}.
\bid{doi={10.1093/biomet/75.2.237}, issn={0006-3444}, mr={0946049}}
\bptok{imsref}%
\end{barticle}
%
\endbibitem

%b27 #&#
\bibitem[\protect\citeauthoryear{Owen}{2001}]{Owen2001}
%
\begin{bbook}[auto:STB|2013/06/05|13:45:01]
\bauthor{\bsnm{Owen},~\bfnm{A.~B.}\binits{A.~B.}}
(\byear{2001}).
\btitle{Empirical Likelihood}.
\bpublisher{Chapman \& Hall/CRC}, \blocation{New York}.
\bptok{imsref}%
\end{bbook}
%
\endbibitem

%b28 #&#
\bibitem[\protect\citeauthoryear{Pan, Xie and Shen}{2010}]{Pan2010Biometrics}
%
\begin{barticle}[mr]
\bauthor{\bsnm{Pan},~\bfnm{Wei}\binits{W.}},
\bauthor{\bsnm{Xie},~\bfnm{Benhuai}\binits{B.}} \AND
\bauthor{\bsnm{Shen},~\bfnm{Xiaotong}\binits{X.}}
(\byear{2010}).
\btitle{Incorporating predictor network in penalized regression with
application to microarray data}.
\bjournal{Biometrics}
\bvolume{66}
\bpages{474--484}.
\bid{doi={10.1111/j.1541-0420.2009.01296.x}, issn={0006-341X}, mr={2758827}}
\bptok{imsref}%
\end{barticle}
%
\endbibitem

%b29 #&#
\bibitem[\protect\citeauthoryear{Petrov}{1995}]{Petrov1995}
%
\begin{bbook}[mr]
\bauthor{\bsnm{Petrov},~\bfnm{Valentin~V.}\binits{V.~V.}}
(\byear{1995}).
\btitle{Limit Theorems of Probability Theory:
Sequences of Independent Random Variables}.
\bseries{Oxford Studies in Probability}
\bvolume{4}.
\bpublisher{Oxford Univ. Press}, \blocation{New York}.
\bid{mr={1353441}}
\bptok{imsref}%
\end{bbook}
%
\endbibitem

%b30 #&#
\bibitem[\protect\citeauthoryear{Qin and Lawless}{1994}]{QinLawless1994AS}
%
\begin{barticle}[mr]
\bauthor{\bsnm{Qin},~\bfnm{Jing}\binits{J.}} \AND
\bauthor{\bsnm{Lawless},~\bfnm{Jerry}\binits{J.}}
(\byear{1994}).
\btitle{Empirical likelihood and general estimating equations}.
\bjournal{Ann. Statist.}
\bvolume{22}
\bpages{300--325}.
\bid{doi={10.1214/aos/1176325370}, issn={0090-5364}, mr={1272085}}
\bptok{imsref}%
\end{barticle}
%
\endbibitem

%b31 #&#
\bibitem[\protect\citeauthoryear{Qu, Lindsay and Li}{2000}]{Quetal2000Bioka}
%
\begin{barticle}[mr]
\bauthor{\bsnm{Qu},~\bfnm{Annie}\binits{A.}},
\bauthor{\bsnm{Lindsay},~\bfnm{Bruce~G.}\binits{B.~G.}} \AND
\bauthor{\bsnm{Li},~\bfnm{Bing}\binits{B.}}
(\byear{2000}).
\btitle{Improving generalised estimating equations using quadratic inference
functions}.
\bjournal{Biometrika}
\bvolume{87}
\bpages{823--836}.
\bid{doi={10.1093/biomet/87.4.823}, issn={0006-3444}, mr={1813977}}
\bptok{imsref}%
\end{barticle}
%
\endbibitem

%b32 #&#
\bibitem[\protect\citeauthoryear{Saulis and
Statulevi{\v{c}}ius}{1991}]{SaulisStatuleviius1991}
%
\begin{bbook}[mr]
\bauthor{\bsnm{Saulis},~\bfnm{L.}\binits{L.}} \AND
\bauthor{\bsnm{Statulevi{\v{c}}ius},~\bfnm{V.~A.}\binits{V.~A.}}
(\byear{1991}).
\btitle{Limit Theorems for Large Deviations}.
\bseries{Mathematics and Its Applications (Soviet Series)}
\bvolume{73}.
\bpublisher{Kluwer Academic}, \blocation{Dordrecht}.
\bnote{Translated and revised from the 1989 Russian original}.
\bid{doi={10.1007/978-94-011-3530-6}, mr={1171883}}
\bptok{imsref}%
\end{bbook}
%
\endbibitem

%b33 #&#
\bibitem[\protect\citeauthoryear{Tang and Leng}{2010}]{TangLeng2010Bioka}
%
\begin{barticle}[mr]
\bauthor{\bsnm{Tang},~\bfnm{Cheng~Yong}\binits{C.~Y.}} \AND
\bauthor{\bsnm{Leng},~\bfnm{Chenlei}\binits{C.}}
(\byear{2010}).
\btitle{Penalized high-dimensional empirical likelihood}.
\bjournal{Biometrika}
\bvolume{97}
\bpages{905--919}.
\bnote{With supplementary material available online}.
\bid{doi={10.1093/biomet/asq057}, issn={0006-3444}, mr={2746160}}
\bptok{imsref}%
\end{barticle}
%
\endbibitem

%b34 #&#
\bibitem[\protect\citeauthoryear{Wang}{2012}]{Wang2012Bioka}
%
\begin{barticle}[mr]
\bauthor{\bsnm{Wang},~\bfnm{H.}\binits{H.}}
(\byear{2012}).
\btitle{Factor profiled sure independence screening}.
\bjournal{Biometrika}
\bvolume{99}
\bpages{15--28}.
\bid{doi={10.1093/biomet/asr074}, issn={0006-3444}, mr={2899660}}
\bptok{imsref}%
\end{barticle}
%
\endbibitem

%b35 #&#
\bibitem[\protect\citeauthoryear{Xue and Zou}{2011}]{XueZou2011bioka}
%
\begin{barticle}[mr]
\bauthor{\bsnm{Xue},~\bfnm{Lingzhou}\binits{L.}} \AND
\bauthor{\bsnm{Zou},~\bfnm{Hui}\binits{H.}}
(\byear{2011}).
\btitle{Sure independence screening and compressed random sensing}.
\bjournal{Biometrika}
\bvolume{98}
\bpages{371--380}.
\bid{doi={10.1093/biomet/asr010}, issn={0006-3444}, mr={2806434}}
\bptok{imsref}%
\end{barticle}
%
\endbibitem

%b36 #&#
\bibitem[\protect\citeauthoryear{Zhao and Li}{2012}]{ZhaoLi2012}
%
\begin{bmisc}[auto:STB|2013/06/05|13:45:01]
\bauthor{\bsnm{Zhao},~\bfnm{S.~D.}\binits{S.~D.}} \AND
\bauthor{\bsnm{Li},~\bfnm{Y.}\binits{Y.}}
(\byear{2012}).
\bhowpublished{Sure screening for estimating equations in ultra-high
dimensions. Unpublished manuscript}.
\bptok{imsref}%
\end{bmisc}
%
\endbibitem

%b37 #&#
\bibitem[\protect\citeauthoryear{Zhu et~al.}{2011}]{Zhuetal2011JASA}
%
\begin{barticle}[mr]
\bauthor{\bsnm{Zhu},~\bfnm{Li-Ping}\binits{L.-P.}},
\bauthor{\bsnm{Li},~\bfnm{Lexin}\binits{L.}},
\bauthor{\bsnm{Li},~\bfnm{Runze}\binits{R.}} \AND
\bauthor{\bsnm{Zhu},~\bfnm{Li-Xing}\binits{L.-X.}}
(\byear{2011}).
\btitle{Model-free feature screening for ultrahigh-dimensional data}.
\bjournal{J. Amer. Statist. Assoc.}
\bvolume{106}
\bpages{1464--1475}.
\bid{doi={10.1198/jasa.2011.tm10563}, issn={0162-1459}, mr={2896849}}
\bptok{imsref}%
\end{barticle}
%
\endbibitem

\end{thebibliography}
\end{document}